\documentclass[10pt]{amsart}
\usepackage{amsmath}
\usepackage{amssymb}
\usepackage{amsfonts}
\usepackage{amsthm}
\usepackage{enumerate}
\usepackage[all]{xy}
\usepackage[latin1]{inputenc}
\usepackage{mathdots}
\usepackage{graphicx}
\usepackage{latexsym}
\usepackage{mathrsfs}
\usepackage{pgf}
\usepackage{braids}
\usepackage{floatrow}
\usepackage{subfig}

\usepackage{tikz}

\PassOptionsToPackage{usenames,dvipsnames,svgnames}{xcolor}

\input xy

\newtheorem{Theo}{Theorem}[section]

\newtheorem{Prop}[Theo]{Proposition}
\newtheorem{Cor}[Theo]{Corollary}
\newtheorem{Lemma}[Theo]{Lemma}

\newtheorem{Conj}[Theo]{Conjecture}

\theoremstyle{definition}

\newtheorem{Exam}[Theo]{Example}

\newtheorem{Defn}[Theo]{Definition}

\def\mystrut(#1,#2){\vrule height #1pt depth #2pt width 0pt}
\newcommand{\rep}{{\rm rep}}

\newcommand{\Hom}{{\rm Hom}}
\newcommand{\Ext}{{\rm Ext}}

\newcommand{\ext}{{\rm ext}}
\newcommand{\SI}{{\rm SI}}

\newcommand{\Z}{\mathbb{Z}}
\newcommand{\C}{\mathcal{C}}
\newcommand{\A}{\mathcal{A}}

\begin{document}

\author{Charles Paquette}
\address{Charles Paquette, Department of Mathematics, University of Connecticut, Storrs, CT, 06269-3009, USA.}
\email{charles.paquette@uconn.edu}
\author{Jerzy Weyman}
\address{Jerzy Weyman, Department of Mathematics, University of Connecticut, Storrs, CT, 06269-3009, USA.}
\email{jerzy.weyman@uconn.edu}

\title[Isotropic Schur roots]{Isotropic Schur roots}

\maketitle

\begin{abstract}
In this paper, we study the isotropic Schur roots of an acyclic quiver $Q$ with $n$ vertices. We study the perpendicular category $\A(d)$ of a dimension vector $d$ and give a complete description of it when $d$ is an isotropic Schur $\delta$. This is done by using exceptional sequences and by defining a subcategory $\mathcal{R}(Q,\delta)$ attached to the pair $(Q,\delta)$. The latter category is always equivalent to the category of representations of a connected acyclic quiver $Q_{\mathcal{R}}$ of tame type, having a unique isotropic Schur root, say $\delta_{\mathcal{R}}$. The understanding of the simple objects in $\A(\delta)$ allows us to get a finite set of generators for the ring of semi-invariants SI$(Q,\delta)$ of $Q$ of dimension vector $\delta$. The relations among these generators come from the representation theory of the category $\mathcal{R}(Q,\delta)$ and from a beautiful description of the cone of dimension vectors of $\A(\delta)$. Indeed, we show that SI$(Q,\delta)$ is isomorphic to the ring of semi-invariants SI$(Q_{\mathcal{R}},\delta_{\mathcal{R}})$ to which we adjoin variables. In particular, using a result of Skowro\'nski and Weyman, the ring SI$(Q,\delta)$ is a polynomial ring or a hypersurface. Finally, we provide an algorithm for finding all isotropic Schur roots of $Q$. This is done by an action of the braid group $B_{n-1}$ on some exceptional sequences. This action admits finitely many orbits, each such orbit corresponding to an isotropic Schur root of a tame full subquiver of $Q$.
\end{abstract}

\section{Introduction}

Let $Q$ be an acyclic quiver with $n$ vertices and $k$ be an algebraically closed field. One crucial tool in representation theory of acyclic quivers is the use of perpendicular categories. If $V$ is a rigid representation of $Q$, then the (left) perpendicular category $^\perp V$ of $V$ is an exact extension-closed abelian subcategory of $\rep(Q)$, where $\rep(Q)$ is the category of finite dimensional representations of $Q$. These perpendicular subcategories were first studied by Geigle and Lenzing in \cite{GL} and also by Schofield in \cite{Sch2}. One important fact is that such a subcategory has a projective generator, or equivalently, it is equivalent to the category of representations of some acyclic quiver. There is a very natural way to generalize perpendicular categories of rigid representations, namely, by taking the perpendicular category $\A(d)$ of a dimension vector $d$; see Section 3. If $V$ is a rigid representation, then $^\perp V = \A(d_V)$ where $d_V$ is the dimension vector of $V$. If $d_V$ is not the dimension vector of a rigid representation, we show that the category $\A(d)$ does not admit a projective generator. However, this category plays a fundamental role for understanding the ring SI$(Q,d)$ of semi-invariants of $Q$ of dimension vector $d$, the latter object being our second object of study in this paper. Indeed, as a ring, SI$(Q,d)$ is generated by the Schofield semi-invariants $C^V$ (see \cite{DW_DetSI}) where $V$ runs through the simple objects in $\A(d)$. In general, $\A(d)$ may have infinitely many non-isomorphic simple objects. Since SI$(Q,d)$ is always finitely generated as a ring, we still need to decide how to pick a nice subset of those generators $C^V$ and find the relations among them. There is no general method, so far, for doing this. When we specialize to the case where $d=\delta$ is an isotropic Schur root of $Q$, then we can answer this problem.
It was proven in \cite{Skowronskiweyman} that when $Q$ is of tame type and $d$ is arbitrary, the ring of semi-invariants ${\rm SI}(Q, d)$ is either a polynomial ring or a hypersurface.
The smallest possible case of hypersurface occurs for the dimension vector $d=\delta$ where $\delta$ is the isotropic Schur root for $Q$ and where $Q$ is of type $\widetilde{\mathbb{E}_6}, \widetilde{\mathbb{E}_7},\widetilde{\mathbb{E}_8}$ or $\widetilde{\mathbb{D}_n}$, for $n \ge 4$. In such a case, the ring of semi-invariants is generated by the $C^V$ where $V$ runs through the quasi-simple exceptional regular representations of $Q$, that is, the exceptional simple objects in $\A(\delta)$. The hypersurface equation in these cases comes from the fact that the subcategory of regular representations in $\rep(Q)$ has exactly three non-homogeneous tubes. If the regular exceptional quasi-simple representations in the first, second and third tubes are respectively denoted $E_1,\ldots ,E_p$, $E'_1,\ldots , E'_q$ and $E''_1,\ldots ,E''_r$, then the hypersurface equation is given by
$$C^{E_1}\cdots C^{E_p}+C^{E'_1}\cdots C^{E'_q}+C^{E''_1}\cdots C^{E''_r}=0.$$
We investigate how much of this structure is preserved for the ring SI$(Q,\delta)$ of an isotropic Schur root $\delta$ of an arbitrary $Q$. We consider the cone of dimension vectors in $\A(\delta)$. Using convex geometry and Radon's theorem, we prove a certain decomposition property of the space of that cone (Proposition\ \ref{PropVectSpacesDecomp}). Using this, we give a complete description of the simple objects in $\A(\delta)$. There exists an exact extension-closed abelian subcategory $\mathcal{R}=\mathcal{R}(Q,\delta)$ of $\rep(Q)$ which has a projective generator, connected and of tame type. In particular, $\mathcal{R}$ is equivalent to $\rep(Q_{\mathcal{R}})$ for a connected quiver of tame type $Q_{\mathcal{R}}$ having a unique isotropic Schur root $\delta_{\mathcal{R}}$, that can also be seen as an isotropic Schur root of $Q$. This subcategory $\mathcal{R}$ is built from the data of an exceptional sequence $(M_{n-2}, \ldots, M_1)$ of length $n-2$ of simple objects in $\A(\delta)$. Up to isomorphisms, the simple objects in $\A(\delta)$ are given by the objects $M_{n-2}, \ldots, M_1$ together with the quasi-simple objects in $\mathcal{R}$. In particular, there are finitely many, up to isomorphism, exceptional simple objects in $\A(\delta)$. Using this, we can get an explicit description of our ring of semi-invariants: SI$(Q,\delta)$ is obtained by adjoining variables to the ring of semi-invariants ${\rm SI}(Q_{\mathcal{R}},\delta_{\mathcal{R}})$. In particular, ${\rm SI}(Q, \delta)$ is still a polynomial ring or a hypersurface. The defining equation, in the hypersurface case, again comes from linear dependance of three products of semi-invariants $C^V$ where $V$ runs through the exceptional quasi-simple objects in $\mathcal{R}(Q,\delta)$. One difference between the general case and the case of quivers of tame type is that the isotropic Schur root $\delta_{\mathcal{R}}$ may differ from $\delta$. Moreover, the root $\delta_{\mathcal{R}}$ needs not lie in the interior of the cone of dimension vectors of $\A(\delta)$. For the last part of this paper, we find an algorithm for finding all isotropic Schur roots of $Q$. We restrict to full exceptional sequences of $\rep(Q)$ that are called exceptional of isotropic type: it is an exceptional sequence of the form $E=(X_1, \ldots, X_{n-1}, X_n)$ where there is an integer $i$ such that the thick subcategory generated by $X_{i}, X_{i+1}$ contains an isotropic Schur root $\delta_E$, called the type of $E$. We explain how the braid group $B_{n-1}$ acts on these sequences to get all the isotropic Schur roots of $Q$. The action admits finitely many orbits, and each orbit contains an exceptional sequence whose type is an isotropic Schur root of a tame full subquiver of $Q$.

\section{Preliminaries}

In this paper, $Q=(Q_0, Q_1)$ is always a connected acyclic quiver with $n$ vertices, unless otherwise indicated, and $k$ is an algebraically closed field. We let $\rep(Q)$ denote the category of ($k$-linear) finite dimensional representations of $Q$ over $k$.

\subsection{Bilinear form, Schur roots} Given a representation $M$, we denote by $d_M$ its dimension vector, which is an element in $(\Z_{\ge 0})^n$. We denote by $\langle -,- \rangle_Q$, or simply by $\langle -, - \rangle$ when there is no risk of confusion, the \emph{Euler-Ringel} form for $Q$. This is the unique $k$-bilinear form in $\mathbb{R}^n$ such that for $M,N \in \rep(Q)$, we have $$\langle d_M, d_N \rangle = {\rm dim}_k \Hom(M,N) - {\rm dim}_k \Ext^1(M,N).$$ Recall that $M \in \rep(Q)$ is called a \emph{Schur representation} if $\Hom(M,M)$ is one dimensional. It is well known from Kac \cite[Lemma 2.1]{Kac0} that in this case, $\langle d_M, d_M \rangle$ is at most one. A Schur representation $M$ is called an \emph{exceptional representation} if $\langle d_M, d_M \rangle = 1$ or, equivalently, $\Ext^1(M,M)=0$.
A dimension vector $d \in (\Z_{\ge 0})^n$ is a \emph{Schur root} if $d = d_M$ for some Schur representation $M$. We call such a $d$ \emph{real, isotropic or imaginary} if $\langle d,d\rangle$ is one, zero or negative, respectively. In this paper, $\delta$ will always denote an isotropic Schur root of $Q$.

\medskip

\subsection{Geometry and semi-invariants} For an arrow $\alpha \in Q_1$, we denote by $t(\alpha)$ its \emph{tail} and by $h(\alpha)$ its \emph{head}. We write ${\rm Mat}_{u\times v}(k)$ for the set of all $u\times v$ matrices over $k$. For a dimension vector $d=(d_1, \ldots, d_n)$, we denote by $\rep(Q,d)$ the space of representations of dimension vector $d$ with fixed vector spaces, that is,
$$\rep(Q,d) = \prod_{\alpha \in Q_1}{\rm Mat}_{t(\alpha)\times h(\alpha)}(k).$$
This space is an affine space and the reductive group GL$_d(k):= \prod_{i=1}^n{\rm GL}_{d_i}(k)$ acts on it by simultaneous conjugation, so that for $M \in \rep(Q,d)$, the GL$_d(k)$-orbit of $M$ is the set of all representations in $\rep(Q,d)$ that are isomorphic to $M$. Since $Q$ is acyclic, the ring of invariants $k[\rep(Q,d)]^{{\rm GL}_d(k)}$ is trivial. Instead, one rather considers the ring of invariants  $k[\rep(Q,d)]^{{\rm SL}_d(k)}$, where $${\rm SL}_d(k):= \prod_{i=1}^n{\rm SL}_{d_i}(k)\subset {\rm GL}_d(k).$$ This ring, denoted SI$(Q,d)$, is called the \emph{ring of semi-invariants} of $Q$ of dimension vector $d$. It is always finitely generated, since ${\rm SL}_{d}(k)$ is reductive. Let $\Gamma$ denote the group of all homomorphisms $\Z^n \to \Z$, which is identified to the set of multiplicative characters of GL$_d(k)$. A well known result states that SI$(Q,d)$ admits a weight space decomposition $${\rm SI}(Q,d) = \oplus_{\tau \in \Gamma}{\rm SI}(Q,d)_\tau$$
where ${\rm SI}(Q,d)_\tau$ is called the space of \emph{semi-invariants of weight }$\tau$. This makes ${\rm SI}(Q,d)$ a $\Gamma$-graded ring.

\medskip

If $V$ is a representation with $\langle d_V, d \rangle = 0$ and $f_V: P_1 \to P_0$ denotes a fixed projective resolution of $V$, then for $M \in \rep(Q,d)$, we have a $k$-linear map $d(V,M) = \Hom(f_V, M)$ given by a square matrix. We define $C^V(-)$ to be the polynomial function on $\rep(Q,d)$ that takes a representation $M$ to the determinant of $d(V,M)$. If we change the projective resolution of $V$, $C^V(-)$ only changes by a non-zero scalar. This function $C^V(-)$ is a semi-invariant of weight $\langle d_V, - \rangle$, and will be called \emph{determinantal semi-invariant}; see \cite{DW_DetSI, Sch2}. The following theorem was proven by Derksen and Weyman in \cite{DW_DetSI} and also by Schofield and van den Bergh in \cite{SchVDB} in the characteristic zero case.

\begin{Theo}[Derksen-Weyman, Schofield-Van den Bergh] \label{ThmDerksenWeyman}Let $d$ be a dimension vector. Then the ring ${\rm SI}(Q,d)$ is spanned over $k$ by the determinantal semi-invariants $C^V(-)$ where $\langle d_V, d \rangle = 0$.
\end{Theo}

There is a dual way to define semi-invariants. One can also take a representation $W$ with $\langle d, d_W \rangle = 0$ and $f_W: P_1' \to P_0'$ a fixed projective resolution of $W$. For $M \in \rep(Q,d)$, we have a $k$-linear map $d(M,W) = \Hom(M, f_W)$ given by a square matrix. We define $C^-(W)$ the polynomial function on $\rep(Q,d)$ that takes a representation $M$ to the determinant of $d(M,W)$. This function $C^-(W)$ is a semi-invariant of weight $-\langle -, d_W\rangle$. As in the above theorem, the ring ${\rm SI}(Q,d)$ is spanned over $k$ by the semi-invariants $C^-(W)$ where $\langle d, d_W \rangle = 0$.

\medskip

\subsection{Exceptional sequences, thick subcategories and braid groups} One of our main tools in this paper will be to make use of particular exceptional sequences. Recall that a sequence
$$(X_1, \ldots, X_r)$$
of exceptional representations is an \emph{exceptional sequence} if $$\Hom(X_i, X_j) = 0 = \Ext^1(X_i, X_j)$$ whenever $i < j$. It is \emph{full} if $r = n$. For such an exceptional sequence $E$, we denote by $\C(E)$ the thick subcategory of $\rep(Q)$ generated by the objects in $E$. Here \emph{thick} means full, closed under extensions, under direct summands, under kernels of epimorphisms, and under cokernels of monomorphisms. The following is well known. We include a proof for the sake of completeness.

\begin{Prop} \label{PropThick}
A full subcategory $\A$ of $\rep(Q)$ is thick if and only if it is exact abelian and extension-closed.
\end{Prop}

\begin{proof}The sufficiency is clear. Suppose that $\mathcal{A}$ is thick. We only need to show that $\mathcal{A}$ has kernels and cokernels.  We will prove this by showing that the kernel and cokernel in $\rep(Q)$ of a morphism in $\A$ lie in $\A$. Let $f: M \to N$ be a morphism in $\mathcal{A}$ and let $u: K \to M$, $v: N \to C$ and $g: M \to E$ be the kernel, cokernel and coimage in $\rep(Q)$. Since $\rep(Q)$ is hereditary and since we have a monomorphism $g': E \cong {\rm im}(f) \to N$, we have a surjective map $\Ext^1(g',K): \Ext^1(N,K) \to \Ext^1(E,K)$.
The short exact sequence $0 \to K \to M \to E \to 0$ is an element in $\Ext^1(E,K)$ and hence is the image of an element in
$\Ext^1(N,K)$.  We have a pullback diagram
$$\xymatrix{0 \ar[r] & K \ar[r]^u \ar@{=}[d] & M \ar[r]^g \ar[d] & E \ar[r] \ar[d]^{g'} & 0 \\
0 \ar[r] & K \ar[r] & M' \ar[r] & N \ar[r] & 0}.$$
This gives rise to a short exact sequence
$$0 \to M \to M' \oplus E \to N \to 0.$$
Since $M,N \in \mathcal{A}$ and $\mathcal{A}$ is thick, we get $E \in \mathcal{A}$.  Hence, $K,C \in \mathcal{A}$.
\end{proof}

Recall that if $E$ is an exceptional sequence of length $r \le n$, then $\C(E)$ is equivalent to the category $\rep(Q_E)$ of representations of an acyclic quiver $Q_E$ with $r$ vertices; see \cite{Sch2}. Denote by $S_1, \ldots, S_r$ the non-isomorphic simple objects in $\C(E)$. The Euler-Ringel form of $\rep(Q)$, restricted to the subgroup of $\Z^n$ generated by the $d_{S_i}$, is isometric to the Euler-Ringel form of $\rep(Q_E)$. In other words, there is an equivalence $\psi: \rep(Q_E) \to \C(E)$ of categories such that for $X,Y \in \rep(Q_E)$, we have $$\langle d_{\psi(X)}, d_{\psi(Y)} \rangle_{Q} = \langle d_X,d_Y \rangle_{Q_E}.$$
In this way, we will often identify a dimension vector for $Q_E$ to a dimension vector for $Q$. Moreover, using the above equivalence, a Schur root for $Q_E$ will be thought of as a Schur root for $Q$.

Let $E = (X_1, \ldots, X_r)$ be an exceptional sequence. If $j > 1$, we denote by $L_{X_{j-1}}(X_j)$ the reflection of $X_j$ to the left of $X_{j-1}$. This is the unique exceptional representation such that we have an exceptional sequence
$$(X_1, \ldots, L_{X_{j-1}}(X_j), X_{j-1}, X_{j+1}, \ldots, X_r).$$
If $j<r$, we denote by $R_{X_{j+1}}(X_j)$ the reflection of $X_j$ to the right of $X_{j+1}$. This is the unique exceptional representation such that we have an exceptional sequence
$$(X_1, \ldots, X_{j-1}, X_{j+1}, R_{X_{j+1}}(X_j), X_{j+2}, \ldots, X_r).$$
These reflections actually have another interpretation in terms of the braid group. Let $r = n$. For $1 \le i \le n-1$, denote by $\sigma_i$ the operation that takes an exceptional sequence $E$ and reflect the $(i+1)$th object to the left of the $i$th one. Hence, $\sigma_1, \ldots \sigma_{n-1}$ act on the set of all exceptional sequences. It is not hard to check that for $1 \le i \le n-2$ and an exceptional sequence $E$, we have $$(\sigma_i (\sigma_{i+1} (\sigma_iE))) = (\sigma_{i+1} (\sigma_{i} (\sigma_{i+1}E)))$$ and for $|i-j|\ge 2$, we have
$$(\sigma_i(\sigma_jE)) = (\sigma_i(\sigma_jE)).$$
Therefore, the braid group
$$B_n := \langle \sigma_1, \ldots, \sigma_{n-1} \mid \sigma_i \sigma_{i+1} \sigma_i = \sigma_{i+1} \sigma_{i} \sigma_{i+1}, \sigma_i\sigma_j = \sigma_j\sigma_i \;\; \text{for}\;|i-j|\ge 2 \rangle$$
on $n$ strands
acts on the set of all full exceptional sequences. It is well known by a result of Crawley-Boevey \cite{CB} that this action is transitive, meaning that all exceptional sequences lie in a single orbit.

\section{Perpendicular subcategories, stable and semi-stable representations}

In this section, we consider perpendicular subcategories of dimension vectors, study them and explain why these categories are related to semi-invariants of quivers. Let $d$ be any dimension vector. Consider
$$\A(d):=\{X \in \rep(Q) \mid \Hom(X,M)=0=\Ext^1(X,M) \; \text{for some} \; M \in \rep(Q,d)\},$$
seen as a full subcategory. This category is called the \emph{(left) perpendicular subcategory} of $d$.

\begin{Prop} \label{PropAd}The subcategory $\A(d)$ is a thick subcategory of $\rep(Q)$, hence exact extension-closed abelian.
\end{Prop}

\begin{proof}
It is clear that $\A(d)$ is closed under direct summands. Let $0 \to X_1 \to X_2 \to X_3 \to 0$ be a short exact sequence in $\rep(Q)$ with exactly two terms $X_s, X_t$ in $\A(d)$. There are $N_s, N_t \in \rep(Q,d)$ such that $\Hom(X_i, N_i)=0 = \Ext^1(X_i, N_i)$ whenever $i=s,t$. Now, for $i=s,t$, there is an open set $\mathcal{U}_i \in \rep(Q,d)$ such that for $Z_i \in \mathcal{U}_i$, we have $\Hom(X_i, Z_i)=0 = \Ext^1(X_i, Z_i)$. Since $\mathcal{U}_s \cap \mathcal{U}_t$ is non-empty, take $N$ a representation in $\mathcal{U}_s \cap \mathcal{U}_t$. Applying $\Hom(-,N)$ to the above exact sequence, we get that the third term lies in $\A(d)$.
\end{proof}

Clearly, $\A$ has a projective generator if and only if it is equivalent to a module category. Since $\A$ is a Hom-finite hereditary abelian category over an algebraically closed field, this means that $\A$ has a projective generator if and only if it is equivalent to $\rep(Q')$ for some finite acyclic quiver $Q'$.

\begin{Lemma} \label{LemmaRank2Perp}Let $d$ be an imaginary or isotropic Schur root in $\rep(Q)$ with $n=2$. Then $\A(d)$ is not empty.
\end{Lemma}

\begin{proof} The result is trivial if $d$ is isotropic. Therefore, assume that the quiver contains at least $m \ge 3$ arrows and $d$ is imaginary. Assume also that the vertices are $\{1,2\}$ with $1$ being the sink vertex. It is sufficient to prove that the result holds for some dimension vector $f$ in the $\tau$-orbit of $d$, where $\tau$ denotes the Coxeter transformation. Consider the cone $C = \{x \in (\Z_{\ge 0})^2 \mid \langle x, x \rangle < 0\}$. The Coxeter transformation is clearly $C$-invariant. Observe that $\tau$ sends the vector $[m-1, 1]$ to $[1, m-1]$. Set $z = [m-1,1]$ and for $i \in \Z$, set $z_i = \tau^iz$. We claim that for $i \in \Z$, the cones $[z_i, z_{i+1}]$, $[z_{i+1}, z_{i+2}]$ only intersect at the ray generated by $z_{i+1}$. Assume otherwise. By continuity of $\tau$ and since $\tau z_i = z_{i+1}$, there exists a ray in $[z_{i}, z_{i+1}]$ that is fixed by $\tau$. Therefore, there is an eigenvector $v$ in $[z_{i}, z_{i+1}]$. Then $\langle v, v \rangle = - \langle v, \tau v \rangle = -\lambda \langle v, v \rangle$ where $\lambda$ is the corresponding eigenvalue. It is well known, see \cite{Pena}, that $\tau$ has a real eigenvalue greater than $1$ and a positive real eigenvalue less than one. Therefore, $\lambda \ne -1$ meaning that $\langle v, v \rangle = 0$ so $v \not \in C$, a contradiction. This proves our claim. This also proves that the two linearly independent eigenvectors of $\tau$ lie on the two boundary rays of $C$. Now, the vectors $z_1, z_2, z_3, \ldots$ converge to one such ray and $z_0, z_{-1}, z_{-2},\ldots$ converge to the other ray. Therefore, we see that the dimension vectors in the half open cone $[z_0, z_1)$ forms a fundamental domain in $C$ for the action of $\tau$. Hence, we may assume that our dimension vector $d$ lies in $[z_0, z_1]$.

Assume first that $d = [q,p]$ with $1 \le p \le q$. Since $d$ lies in $[z_0, z_1]$, we have $1 \le \frac{q}{p} \le m-1$. By \cite[Theo. 1]{BuDr}, the Hilbert null-cone of $\rep(Q,d)$ is not the entire space if $m > \lceil\frac{q}{p}\rceil$ and $p > 1$. The first condition always holds since $d \in [z_0, z_1]$. So if $p > 1$, then there has to be a representation $M$ of dimension vector $d$ and a semi-invariant that does not vanish at $M$. By Theorem \ref{ThmDerksenWeyman}, this yields a representation $N$ with $C^N(M)\ne 0$, meaning that $N \in \A(d)$. Assume now that $p=1$ and take $M$ a general representation of dimension vector $d$. Consider a general representation $Z$ of dimension vector $[mq-1, q]$. By construction $\langle d_Z, d_M \rangle = 0$. Now, a proper subrepresentation of $M$ has dimension vector $d_i=[i,0]$ for $0 \le i \le q$. Notice that $\langle d_Z, d_i \rangle \le 0$. Therefore, it follows from King's criterion \cite{King} that $M$ is $\langle d_Z, - \rangle$-semistable. Therefore, there is a positive integer $r$ and a representation $Z'$ of dimension vector $rd_Z$ such that $C^{Z'}(M) \ne 0$, meaning that $Z' \in \A(d)$. Similarly, one can prove that $\A(d)$ contains a non-zero representation if $d = [q,p]$ with $1 \le q \le p$.
\end{proof}

For two dimension vectors $d_1,d_2$, let $\ext(d_1,d_2)$ denote the minimal value of ${\rm dim}_k \Ext^1(M_1,M_2)$ where $(M_1,M_2) \in \rep(Q,d_1) \times \rep(Q,d_2)$. Similarly, we let $\hom(d_1,d_2)$ denote the minimal value of ${\rm dim}_k \Hom(M_1,M_2)$ where $(M_1,M_2) \in \rep(Q,d_1) \times \rep(Q,d_2)$. There is an open subset $\mathcal{U}_1$ of $\rep(Q,d_1)$ and an open subset $\mathcal{U}_2$ of $\rep(Q,d_2)$ such that for $M_1 \in \mathcal{U}_1, M_2 \in \mathcal{U}_2$, we have
$$\langle d_1, d_2 \rangle = {\rm dim}_k \Hom(M_1,M_2) -{\rm dim}_k \Ext^1(M_1,M_2) = \hom(d_1,d_2) - \ext(d_1,d_2).$$
We write $d_1 \bot d_2$ if $\hom(d_1,d_2) = \ext(d_1,d_2) = 0$. Observe that $d_1 \bot d_2$ implies $\langle d_1, d_2 \rangle = 0$. A sequence $(d_1, \ldots, d_r)$ of Schur roots with $d_i \bot d_j$ whenever $i < j$ is called an \emph{orthogonal sequence of Schur roots}.

Let $d$ be a dimension vector. Due to results of Kac \cite{Kac}, there is a decomposition, denoted $d=\alpha_1 \oplus \cdots \oplus \alpha_m$, having the property that there exists an open (dense) subset $\mathcal{U}_d$ of $\rep(Q,d)$ such that for $M \in \mathcal{U}_d$, we have a decomposition $M \cong M_1 \oplus \cdots \oplus M_m$, where each $M_{i}$ is a Schur representation with $d_{M_{i}} = \alpha_i$. Moreover, $\Ext^1(M_i,M_j)=0$ when $i \ne j$. The latter decomposition of $d$ is unique up to ordering, and is called the \emph{canonical decomposition} of $d$. The dimension vectors $\alpha_i$ are clearly Schur roots, however, they do not need be distinct. Sometimes, it is more convenient to write the above decomposition as
$$(*) \quad d=p_1\beta_1 \oplus \cdots \oplus p_r\beta_r,$$
where the $\beta_i$ are pairwise distinct and $p_i$ is the number of $1 \le j \le m$ with $\beta_i = \alpha_j$. It follows from \cite[Theorem 3.8]{SchofieldGeneral} that when $\beta_i$ is imaginary, then $p_i=1$. When writing a canonical decomposition as in $(*)$, we adopt the convention that when $\alpha$ is an imaginary Schur root and $p$ is a positive integer, then $p\alpha$ is just one root (not $p$ times the root $\alpha$ as when $\alpha$ is real or isotropic). With this convention, Schofield has proven in \cite[Theorem 3.8]{SchofieldGeneral} the following result.

\begin{Prop}[Schofield]\label{CanDecompSchofield} Let $d=p_1\beta_1 \oplus \cdots \oplus p_r\beta_r$ be the canonical decomposition of $d$. If $p$ is a positive integer, then $pd=pp_1\beta_1 \oplus \cdots \oplus pp_r\beta_r$ is the canonical decomposition of $pd$, using the above convention for imaginary Schur roots.
\end{Prop}

The following generalizes Lemma \ref{LemmaRank2Perp}.

\begin{Lemma} \label{LemmaPerp} Let $d$ be be a dimension vector in $\rep(Q)$ whose canonical decomposition involves an isotropic Schur root or an imaginary Schur root. Then $\A(d)$ is not empty.
\end{Lemma}

\begin{proof} We follow the algorithm in \cite{DWAlgo} to find the canonical decomposition of $d$. In particular, we apply this algorithm until the first step where an imaginary or isotropic Schur root is created. In particular, there is an orthogonal sequence $(\alpha_1, \ldots, \alpha_r)$ of real Schur roots with positive integers $p_1, \ldots, p_r$ such that $d = \sum_{i=1}^rp_i\alpha_i$. Moreover, there is a pair $\alpha_t, \alpha_{t+1}$ with $\langle \alpha_{t+1}, \alpha_t \rangle < 0$ and $\gamma = p_t\alpha_t + p_{t+1}\alpha_{t+1}$ is such that $\langle \gamma, \gamma \rangle \le 0$.

Assume first that $r<n$. Then we can find a real Schur root $\alpha_0$ such that $\alpha_0 \bot \alpha_j$ for $1 \le j \le r$. Then, an exceptional representation of dimension vector $\alpha_0$ lies in $\A(d)$. Therefore, we may assume that $r = n$. Let $\gamma' = \gamma$ if $\langle \gamma, \gamma \rangle < 0$ and $\gamma'$ be the smallest indivisible dimension vector in the ray of $\gamma$, otherwise. Set $p \ge 1$ with $\gamma = p \gamma'$. Observe that $\gamma'$ is a Schur root. The next step of the algorithm replaces the sequence $(\alpha_1, \ldots, \alpha_n)$ with positive integers $p_i$ by the orthogonal sequence $(\alpha_1, \ldots, \alpha_{t-1}, \gamma', \alpha_{t+2},\ldots,\alpha_n)$ of Schur roots with positive integers $p_1, \ldots, p_{t-1},p,p_{t+2},\ldots, p_n$. Now, $$d = \sum_{i=1}^{t-1}p_i\alpha_i + p\gamma' + \sum_{i=t+2}^np_i\alpha_i.$$ For $1 \le i \le n$, let $M_i$ denote an exceptional representation of dimension vector $\alpha_i$. The root $\gamma'$ is a root in $\C(M_t, M_{t+1})$. Since $\C(M_t, M_{t+1})$ is equivalent to the category of representations of a quiver with two vertices, it follows from Lemma \ref{LemmaRank2Perp} that there is a dimension vector $\nu$ in $\C(M_t, M_{t+1})$ with $\nu \bot \gamma'$. Clearly, $\langle \nu, \nu \rangle < 0$ and hence $\nu$ is an imaginary Schur root. Let $Z$ be a general representation of dimension vector $\nu$.
Since $Z$ is in general position, we may assume that it lies in $M_j^\perp$ for $1 \le j \le t-1$ and in $^\perp M_j$ for $t+2 \le j \le n$. Therefore, we may assume it lies in $\C(M_t, M_{t+1})$. Since $Z,M_{t-1}$ are in general position and $\Ext^1(M_{t-1},Z)=0$, either $\Hom(Z,M_{t-1})=0$ or $\Ext^1(Z,M_{t-1})=0$ by \cite[Theo. 4.1]{SchofieldGeneral}. We construct a representation $Z'$ that lies in $M_j^\perp$ for $1 \le j \le t-2$, in $^\perp M_j$ for $j = t-1$ and $t+2 \le j \le n$ and lies in $^\perp N$ where $N$ is a general representation of dimension vector $\gamma'$.

Assume first that $\Hom(Z,M_{t-1})=0$. It follows from \cite[Cor. 14]{DWAlgo} that $\nu'=\nu - \langle \nu, \alpha_{t-1}\rangle \alpha_{t-1}$ is a Schur root with $\hom(\nu', \alpha_{t-1})=0$. Moreover, since $\langle \nu', \alpha_{t-1}\rangle = 0$, we also have $\ext^1(\nu', \alpha_{t-1})=0$. Therefore, we have an orthogonal sequence $$(\alpha_1, \ldots, \alpha_{t-2}, \nu', \alpha_{t-1}, \gamma', \alpha_{t+2}, \ldots, \alpha_n)$$ of Schur roots. Consider $Z'$ a general representation of dimension vector $\nu'$. Then $Z'$ satisfies the above wanted conditions for $Z'$.

Assume now that $\Ext^1(Z,M_{t-1})=0$. It follows from Lemma 2.3 in \cite{SchofieldGeneral} that the non-zero morphisms $Z \to M_{t-1}$ are either all injective or all surjective. Let $f = {\rm dim}_k \Hom(Z,M_{t-1})$. Assume first that all non-zero morphisms $Z \to M_{t-1}$ are surjective. We get an epimorphism $Z \to M_{t-1}^f$ given by the basis elements of $\Hom(Z,M_{t-1})$. It is not hard to check that the kernel $Z'$ lies in $^\perp M_{t-1}$. Clearly, $Z'$ satisfies the above wanted conditions for $Z'$. Similarly, we can construct such a $Z'$ if the non-zero morphisms $Z \to M_{t-1}$ are injective. We can continue this process and construct a representation $X$ that lies in $^\perp M_j$ for $1 \le j \le t-1$ and $t+2 \le j \le n$ and lies in $^\perp N$ where $N$ is a general representation of dimension vector $\gamma'$. In particular, $X$ lies in $\A(d)$. Observe that $X$ is not the zero representation since $\alpha_1, \ldots, \alpha_{t-1}, \nu, \gamma'$ are linearly independent.
\end{proof}

Recall that a representation $V \in \rep(Q)$ is \emph{rigid} if $\Ext^1(V,V)=0$. Hence, the indecomposable rigid representations are the exceptional ones. Observe that there is a one-to-one correspondence
$$\{\text{real Schur roots}\} \leftrightarrow \{\text{iso. classes of exceptional representations}\}.$$
The dimension vector of a rigid representation is called \emph{prehomogeneous}. Observe that if $V$ is rigid, then the GL$(d_V)$-orbit of $V$ is open in $\rep(Q,d_V)$. In this case, the canonical decomposition of $d$ involves the dimension vectors of its indecomposable direct summands, so involves only real Schur roots. Conversely, if the canonical decomposition of a dimension vector $d$ involves only real Schur roots, then $\rep(Q,d)$ has an open orbit and hence, $d = d_V$ where $V$ is rigid. Thus, the above correspondence extends to the following one-to-one correspondence:
$$\{d \mid d \; \text{is prehomogeneous}\} \leftrightarrow \{\text{iso. classes of rigid representations}\}.$$

\begin{Prop} The subcategory $\A(d)$ has a projective generator if and only if $d$ is prehomogeneous. In this case,
$$\A(d) = ^\perp V = \{X \in \rep(Q) \mid \Hom(X,V)=0=\Ext^1(X,V)\}$$
where $V$ is rigid with $d = d_V$.
\end{Prop}

\begin{proof}
Assume that $\A(d)$ has a projective generator. Therefore, it is equivalent to the category of representations of an acyclic quiver. Let $M_1, \ldots, M_r$ be the indecomposable simple objects in $\A(d)$, up to isomorphism. We may assume that they are ordered so that $E:=(M_1, \ldots, M_r)$ is an exceptional sequence in $\A(d)$ and hence, also an exceptional sequence in $\rep(Q)$. We have $\A(d) = \C(E)$. As argued in the proof of Proposition \ref{PropAd}, for $1 \le i \le r$, there are open sets $\mathcal{U}_i$ in $\rep(Q,d)$ such that for $N_i \in \mathcal{U}_i$, we have $\Hom(M_i,N_i)=0=\Ext^1(M_i,N_i)$. Now, $\bigcap\mathcal{U}_i$ is non-empty and we let $N$ lies in it. Then, for $1 \le i \le r$, we have $\Hom(M_i,N)=0=\Ext^1(M_i,N)$. Since the $M_i$ are the simple objects in $\C$, it follows that $N \in \C(E)^\perp$. Thus, $\A(d) \subseteq ^\perp N$. However, by definition of $\A(d)$, we have $^\perp N \subseteq \A(d)$. Therefore, $\A(d) = ^\perp N$. Observe that $\C(E)^\perp$ is equivalent to the category of representations of an acyclic quiver $Q'$. We can think of $N$ as a representation in $\rep(Q')$ with dimension vector $d$. Assume that $N$ is not rigid. Then the canonical decomposition of $d$ (as a dimension vector of $Q'$) involves an isotropic or imaginary Schur root of $Q'$. It follows from Lemma \ref{LemmaPerp} that there is a representation $Z$ in the category $\A(d)$ for $\rep(Q')$. This means that $Z \in \A(d)$ but $Z \not \in \C(E)$, a contradiction. Therefore, $N$ is rigid and this proves the necessity. Assume now that $d = d_V$ where $V$ is rigid. Since $V$ is in general position, any $Z \in \A(d)$ has to be in $^\perp V$ and hence, $\A(d) = ^\perp V.$
\end{proof}

Given a dimension vector $d$, we fix $\sigma_d:= -\langle -, d \rangle$, which is called the \emph{weight associated} to $d$.
Recall that $M \in \rep(Q)$ is $\sigma_d$-\emph{semistable} if there is a positive integer $m$ and a semi-invariant $f$ of weight $m\sigma_d$ in SI$(Q, d_M)$ such that $f$ does not vanish at $M$. If $M$ is $\sigma_d$-semistable and has no proper (non-zero) $\sigma_d$-semistable subobject, then it is $\sigma_d$-\emph{stable}. It follows from King's criterion \cite{King} that $M$ is $\sigma_d$-semistable (resp. $\sigma_d$-stable) if and only if $\sigma_d(d_M) = 0$ and $\sigma_d(f) \le  0$ (resp. $\sigma_d(f) <  0$) whenever $M$ has a proper non-zero subobject of dimension vector $f$. We will see that this notion of semistability is related to perpendicular subcategories. Observe first that for any dimension vector $d$, and a positive integer $m$, we have $\A(d) \subseteq \A(md)$. The other inclusion is not true, in general, if the canonical decomposition of $d$ involves an imaginary Schur root.

\begin{Prop}
Let $d$ be a dimension vector. The (simple) objects in $\cup_{m \ge 1}\A(md)$ are the $\sigma_d$-(semi)stable objects. If the canonical decomposition of $d$ does not involve imaginary Schur roots, then $\cup_{m \ge 1}\A(md) = \A(d)$.
\end{Prop}

\begin{proof} Let $M \in \A(d)$. Then there exists $M(d) \in \rep(Q,d)$ such that $$\Hom(M,M(d)) = 0 = \Ext^1(M,M(d)).$$ Therefore, the semi-invariant $C^-(M(d))$ of weight $\sigma_d$ in SI$(Q, d_M)$ does not vanish at $M$. Since there exists a semi-invariant in ${\rm SI}(Q,d_M)_{\sigma_d}$ that does not vanish on $M$, $M$ is $\sigma_d$-semistable. Conversely, assume that $M$ is $\sigma_d$-semistable. Then there exists $m \ge 1$ and a semi-invariant $f\in{\rm SI}(Q,d_M)_{m \sigma_d}$ that does not vanish on $M$. Now, $f$ is given by a semi-invariant of the form $C^-(N)$ for some representation $N$ of dimension vector $m d$. We see that $M$ lies in $\A(md) \subseteq \cup_{i \ge 1}\A(id)$. For $\sigma_d$-stable representations, one just needs to use king's criterion, as mentioned in the paragraph before this proposition.

Assume now that the canonical decomposition of $d$ does not involve imaginary Schur root. It follows from Proposition \ref{CanDecompSchofield} that a general representation of dimension vector $md$ has a direct summand of dimension vector $d$. Let $M \in \A(md)$. Then there exists $N \in \rep(Q,md)$ such that $C^M(N) \ne 0$. Also, $N$ can be taken to be in general position. By the above observation, there is a summand $N'$ of $N$ of dimension vector $d$. Thus $C^M(N')\ne 0$, meaning that $M \in \A(d)$. Therefore, $\A(md) = \A(d)$ for all $m \ge 1$.
\end{proof}

The following proposition explains why the study of perpendicular subcategories is directly related to the study of semi-invariants.

\begin{Prop}
The ring $\SI(Q,d)$ is generated by the semi-invariants $C^V(-)$ where $V$ is simple in $\A(d)$.
\end{Prop}

\begin{proof}
We know from Theorem \ref{ThmDerksenWeyman} that the ring $\SI(Q,d)$ is generated, over $k$, by semi-invariants of the form $C^V(-)$ where $V$ is a representation with $\langle d_V, d \rangle = 0$. If $\Hom(V,M) \ne 0$ for all $M \in \rep(Q, d)$, then $C^V(-)$ is the zero semi-invariant. Otherwise, $V$ lies in $\A(d)$. If $V$ is not a simple object in $\A(d)$, then there exists a simple subobject $V_1$ of $V$ in $\A(d)$. Since $\A(d)$ is thick, this yields a short exact sequence
$$0 \to V_1 \to V \to V_2 \to 0$$
in $\A(d)$ and it follows from \cite[Lemma 1]{DW_DetSI} that $C^{V}(-)=aC^{V_1}(-)C^{V_2}(-)$ where $a \in k$ is non-zero. Repeating this reduction for the object $V_2$ yields that $C^V(-)$ is, up to a scalar, a product of semi-invariants as in the statement.
\end{proof}

If $\sigma$ is a weight, a dimension vector $d$ is called \emph{$\sigma$-semistable} (resp. $\sigma$-\emph{stable}) if a general representation of dimension vector $d$ is $\sigma$-semistable (resp. $\sigma$-stable). Given a $\sigma$-semistable dimension vector $d$, there exists a $\sigma$-stable decomposition
$d = d_1 + \cdots + d_r$
of $d$ where a general representation $M$ of dimension vector $d$ has a filtration
$$0 \subset M_1 \subset \cdots \subset M_r = M$$
such that for $1 \le i \le r$, $M_i/M_{i-1}$ has dimension vector $d_i$ and is $\sigma$-stable; see \cite{DWSchurSeqn}. In this case, $d_i$ will be called a \emph{$\sigma$-stable factor} of $d$. The $\sigma$-stable decomposition of a $\sigma$-semistable dimension vector is unique up to ordering. Moreover, it was shown in \cite{DWSchurSeqn} that the pairwise distinct $\sigma$-stable factors $f_1, \ldots, f_s$ of $d$ can be ordered in such a way that $f_i \bot f_j$ whenever $i < j$. So the sequence $(f_1, \ldots, f_s)$ is an orthogonal sequence of Schur roots.

\begin{Lemma} \label{Lemma0} Let $\sigma$ be a weight and $d$ be a Schur root that is $\sigma$-semi-stable. If $d$ is real, then all the $\sigma$-stable factors of $d$ are real Schur roots. If $d$ is isotropic, then all the $\sigma$-stable factors of $d$ are real or isotropic Schur roots.
\end{Lemma}

\begin{proof}
Let $d_1, \ldots, d_r$ be the distinct $\sigma$-stable factors of $d$. We may assume that $(d_1, \ldots, d_r)$ is an orthogonal sequence of Schur roots.
We can use the algorithm in \cite{DWAlgo} for finding the canonical decomposition of $d$ starting with the sequence $(d_1, \ldots, d_r)$.
If one of $d_i$ is imaginary (or isotropic), it follows from the algorithm that $d$ will be imaginary (resp. isotropic or imaginary). Therefore, if $d$ is real, then all $d_i$ are real. If $d$ is isotropic, then no $d_i$ is imaginary.
\end{proof}

Let $E=(X_1, \ldots, X_r)$ be an exceptional sequence and consider the subcategory $\C(E)=\C(X_1, \ldots, X_r)$. Assume that the dimension vector $d$ lies in $\C(E)$, that is, there is an object in $\C(E)$ having $d$ as dimension vector. The stability condition $\sigma_d = -\langle - , d \rangle$ also gives rise to a stability condition, denoted $\sigma_{\C(E), d}$, in $\C(E)$. An object $X \in \C(E)$ is said to be \emph{relative $\sigma_d$-semistable} (resp. \emph{relative $\sigma_d$-stable}) in $\C(E)$ provided it is $\sigma_{\C(E), d}$-semistable (resp. $\sigma_{\C(E),d}$-stable).

\begin{Prop} \label{PropRelativeStable}Let $E=(X_1, \ldots, X_r)$ be an exceptional sequence and let $d$ be a dimension vector lying in $\C(E)$. Let $X \in \C(E)$.
\begin{enumerate}[$(1)$]
    \item The object $X$ is relative $\sigma_{d}$-semistable in $\C(E)$ if and only if it is $\sigma_d$-semistable.
\item If $X$ is $\sigma_d$-stable, then $X$ is relative $\sigma_{d}$-stable in $\C(E)$. \end{enumerate}
\end{Prop}

\begin{proof}
Observe that $X$ is $\sigma_d$-semistable if and only if $-\langle d_X, d \rangle = 0$ and for any subobject $X'$ of $X$, we have $-\langle d_{X'}, d \rangle \le 0$. Similarly, $X$ is $\sigma_{\C(E),d}$-semistable if and only if $-\langle d_X, d \rangle = 0$ and for any subobject $X'$ of $X$ in $\C(E)$, we have $-\langle d_{X'}, d \rangle \le 0$. Now, a subobject of $X$ in $\C(E)$ is also a subobject of $X$ in $\rep(Q)$. This shows the sufficiency of (1). Suppose now that $X$ is $\sigma_{\C(E),d}$-semistable. Then there exists $M$ in $\C(E)$ of dimension vector $md$ for some $m > 0$ such that $C^X(M)\ne 0$, that is, $\Hom(X,M)=\Ext^1(X,M)=0$. Since $M$ has dimension vector $md$ in $\rep(Q)$, this gives that $X$ is $\sigma$-semistable.
Assume now that $X$ is $\sigma$-stable. Then $X$ is $\sigma_{\C,d}$-semistable. If $X$ is not $\sigma_{\C,d}$-stable, then there exists a proper subobject $X'$ of $X$ in $\C(E)$ such that $-\langle d_{X'}, d \rangle = 0$. Now, $X'$ is also a subobject of $X$ in $\rep(Q)$ with $-\langle d_{X'}, d \rangle = 0$, and this contradicts the fact that $X$ is $\sigma_d$-stable.
\end{proof}

\section{The case of an isotropic Schur root}

In this section, $\delta$ stands for an isotropic Schur root. The weight $\sigma_\delta$ will simply be denoted $\sigma$, when there is no risk of confusion. Our aim is to describe all simple objects in $\A(\delta)$ or, equivalently, all $\sigma_\delta$-stable objects. We start with the following proposition; see \cite{Pa}.

\begin{Prop} \label{ExceptionalSequenceIsotropic}
There exists an exceptional sequence $(V,W)$ in $\rep(Q)$ such that $\C(V,W)$ is tame and $\delta = d_V + d_W$ with $\langle d_W, d_V \rangle = -2$. In particular, $^\perp V \cap ^\perp W \subseteq \A(\delta)$.
\end{Prop}

\begin{proof}
As shown in \cite{Pa}, there is an exceptional sequence $E=(V,W)$ of length two such that $\delta$ is a root in $\C(E)$. Since $E$ has length two, $\C(E)$ is equivalent to the category of representations of an acyclic quiver $Q_E$ with two vertices. But $Q_E$ has to have an isotropic Schur root. Therefore, $Q_E$ is the Kronecker quiver. With no loss of generality, we may assume that $V,W$ are the simple objects of $\C(E)$. Therefore, $\langle d_W, d_V \rangle = -{\rm dim}_k\Ext^1(V,W)$ and ${\rm dim}_k\Ext^1(V,W)$ is the number of arrows in $Q_E$. The second part of the statement is trivial.
\end{proof}

\begin{Lemma} \label{Lemma1}
Let $Q$ have at least three vertices. Then there is at least one exceptional simple object in $\A(\delta)$.
\end{Lemma}

\begin{proof}
Suppose that all the $\sigma$-stable dimension vectors are isotropic or imaginary. By Proposition \ref{ExceptionalSequenceIsotropic}, we have an exceptional sequence $(V,W)$, such that $\delta = d_V + d_W$. We can extend this to an exceptional sequence $(U,V,W)$ and hence, $U \in \A(\delta)$. This means that $d_U$ is $\sigma$-semistable. Now, we apply Lemma \ref{Lemma0}.
\end{proof}

The subcategories of the form $\C(A,B,C)$ where $(A,B,C)$ is an exceptional sequence will play a crucial role in our investigation. The following lemma, which is easy to check, provides a description of the quivers with three vertices having an isotropic Schur root.

\begin{Lemma} \label{rank3} Let $(A,B,C)$ be an exceptional sequence such that $\C(A,B,C)$ contains $\delta$. Then either
\begin{enumerate}[$(1)$]
    \item the category $\C(A,B,C)$ is wild and connected or
    \item the category $\C(A,B,C)$ is equivalent to $\rep(Q')$, where $Q'$ is either of type $\widetilde{\mathbb{A}}_{2,1}$ or a union of the Kronecker quiver and a single vertex.
\end{enumerate}
\end{Lemma}

Let us denote by $\tau$ the Auslander-Reiten translation in $\rep(Q)$. An indecomposable representation that lies in the $\tau$-orbit of a projective (resp. injective) representation is called \emph{preprojective} (resp. preinjective). Baer and Strau${\ss}$ have proven the following crucial result; see \cite{Baer} or \cite[Theorem B]{Strauss}.

\begin{Lemma}[Baer, Strau${\ss}$] \label{ResultBaerStrauss}Let $Q$ be of wild type and let $X$ be exceptional. If $X^\perp$ is of finite or tame type, then $X$ has to be preprojective or preinjective.
\end{Lemma}

The following result describes a way to produce other isotropic Schur roots starting with an exceptional sequence $(U,V,W)$ where $\delta = d_V + d_W$.

\begin{Lemma} \label{LemmaIsoPreinj} Let $E=(U,V,W)$ be an exceptional sequence such that $\C(V,W)$ is tame with isotropic Schur root $\delta = d_V + d_W$. Reflect both $V,W$ to the left of $U$ to get an exceptional sequence $(V', W', U)$. Let $\delta'$ be the unique isotropic Schur root in $\C(V', W')$. Then $\delta' = \delta - \langle \delta, d_U \rangle d_U$.
\begin{enumerate}[$(1)$]
    \item If $U$ is preinjective in $\C(E)$, then $Q$ is wild connected and $\langle \delta, d_U \rangle \ge 0$.
    \item If $U$ is preprojective in $\C(E)$, then $Q$ is wild connected and $\langle \delta, d_U \rangle \le 0$.
    \item Otherwise, $\C(E)$ is tame, $U$ is regular or is simple projective-injective in $\C(E)$, $\langle \delta, d_U \rangle = 0$ and $\delta' = \delta$.
\end{enumerate}
\end{Lemma}

\begin{proof} We may assume that $\rep(Q) = \C(E)$. If $Q$ is wild, then it follows from Lemma \ref{rank3} that $Q$ is connected. In this case, by Lemma \ref{ResultBaerStrauss}, $U$ cannot be regular. Assume first that $Q$ is tame. If $Q$ is tame connected, $Q$ is of type $\widetilde{\mathbb{A}}_{2,1}$ and $U$ has to be isomorphic to one of the two quasi-simple regular exceptional representations. If $Q$ is tame disconnected, then we see that $U$ is the simple representation corresponding to the connected component of $Q$ of type $\mathbb{A}_1$. So if $Q$ is tame, it is clear that $\langle \delta, d_U \rangle = 0$ and $\delta = \delta'$. So we may assume that $Q$ is wild connected.
We have an orthogonal sequence of Schur roots $(d_U, \delta)$. Set $\delta'' = \delta - \langle \delta, d_U \rangle d_U$. One easily checks that $\langle \delta'', d_U \rangle = 0$ and $\langle \delta'', \delta'' \rangle = 0$. Since $^\perp U = \C(V', W')$ is of tame type, the Euler-Ringel form $\langle - , - \rangle_{\C(V', W')}$ restricted to $\C(V', W')$ is positive semi-definite. Therefore, $\delta''$ is an integral multiple of the isotropic root $\delta'$ in $\C(V', W')$. Since each of $\delta', \delta''$ is a sum or difference of $d_{V'}, d_{W'}$, we see that $\delta'' = \pm \delta'$. Assume first that $U$ is preinjective. A general representation $M(\delta)$ of dimension vector $\delta$ is regular while a general representation of dimension vector $d_U$ is isomorphic to $U$ hence preinjective. We have $\Ext^1(M(\delta), U)=0$ and hence $\langle \delta, d_U \rangle = {\rm dim}_k\Hom(M(\delta), U) \ge 0$. If $U$ is preprojective, then $\Hom(M(\delta), U) = 0$. This gives $\langle \delta, d_U \rangle \le 0$. Suppose that $\delta' = \langle \delta, d_U \rangle d_U - \delta$. This is only possible if $\langle \delta, d_U \rangle > 0$ and hence, if $U$ is preinjective. Then $\langle \delta, d_U \rangle d_U = \delta + \delta'$. The region $\mathscr{R}$ given by $\langle d, d \rangle \le 0$ is the cone over a two dimensional ellipse and hence, is a convex cone.
Since both $\delta, \delta'$ lie on the boundary of $\mathscr{R}$, we see that $d_U$ lies in $\mathscr{R}$. Hence, $\langle d_U, d_U \rangle \le 0$, a contradiction to $U$ being exceptional. Therefore, $\delta'' = \delta'$ is the wanted isotropic Schur root.
\end{proof}

The root $\delta'$ defined above will be denoted $L_U(\delta)$, as it is the reflection of $\delta$ to the left of the real Schur root $d_U$.

\begin{Prop} \label{Prop3Vertices}
Let $Q$ be wild with three vertices, and let $(V,W)$ be an exceptional sequence of tame type containing $\delta$. Complete this to a full exceptional sequence $(U,V,W)$. Then $U$ is preprojective or preinjective.
\begin{enumerate}[$(1)$]
    \item If $U$ is preprojective, then the $\sigma$-stable representations are, up to isomorphism, $U$ or $M(\delta)$ where $M(\delta)$ is a general representation of dimension vector $\delta$ in $\C(V,W)$. In particular, the $\sigma$-stable dimension vectors are $d_U, \delta$.
    \item If $U$ is preinjective, then the $\sigma$-stable representations are, up to isomorphism, $U$ or $M(\delta')$ where $M(\delta')$ is a general representation of dimension vector $\delta' = L_U(\delta)$ in $\C(V',W')$, where $(V', W', U)$ is exceptional. In particular, the $\sigma$-stable dimension vectors are $d_U, \delta'$.
\end{enumerate}
In both cases, the $\sigma$-stable dimension vectors are the two extremal rays in the cone of $\sigma$-semistable dimension vectors.
\end{Prop}

\begin{proof} Recall that the region $\langle d, d \rangle\le 0$ is the cone over a two dimensional ellipse. The ray $[\delta]$ lies on the quadric $\langle d, d \rangle=0$. Since $\C(V, W)$ is tame, the Euler-Ringel form for $\C(V, W)$ is positive semi-definite. In particular, any linear combination of $d_V, d_W$ that lie on the region $\langle d, d \rangle \le 0$ has to be a multiple of $\delta$. Therefore, the cone $C_{V, W}$ generated by $d_V, d_W$ is the tangent plane of $\langle d, d \rangle=0$ at $\delta$. If $U'$ is exceptional and $\sigma$-semistable, then we have an exceptional sequence $(U', X, Y)$ such that $[\delta]$ lies in the cone $C_{X, Y}$ generated by $d_{X}, d_{Y}$. Therefore, as argued above, $C_{X, Y}$ is the tangent plane of $\langle d, d \rangle=0$ at $\delta$. Hence, $C_{V, W} = C_{X, Y}$ and thus, since $\langle - , - \rangle$ is non-degenerate, the rays of $d_U$, $d_{U'}$ coincide. This gives $U \cong U'$ since $U, U'$ are exceptional. This proves that there is a unique $\sigma$-semistable real Schur root $d_U$ and a unique $\sigma$-semistable exceptional representation, up to isomorphism. Hence, $d_U$ has to be $\sigma$-stable by Lemma \ref{Lemma1}. Observe that $Q$ is connected by Lemma \ref{rank3}.
By Lemma \ref{ResultBaerStrauss}, since $\C(V,W)$ is tame, we know that $U$ is preinjective or preprojective. Let $M$ be a $\sigma$-stable representation that is not isomorphic to $U$.  By the previous argument, $M$ cannot be exceptional. Suppose first that $U$ is preprojective. If $M$ is not in $U^\perp$, then $\Ext^1(U,M)\ne 0$ as $\Hom(U,M)=0$. Therefore, $M$ has to be preprojective, contradicting that $M$ is not exceptional. Hence, $M$ is $\sigma$-stable and lies in $U^\perp = \C(V,W)$. By Proposition \ref{PropRelativeStable}, $M$ is relative $\sigma$-stable in $U^\perp$. Therefore, $M \cong M(\delta)$ is a general representation of dimension vector $\delta$ in $\C(V,W)$. Suppose finally that $U$ is preinjective. If $M$ is not in $^\perp U$, then $\Ext^1(M,U)\ne 0$ and this means that $M$ is preinjective, thus exceptional, a contradiction. Therefore, $M$ lies in $^\perp U$. Now, $^\perp U$ is also tame. Consider the exceptional sequence $(V', W', U)$ where $V', W'$ are obtained from $V,W$ by reflecting to the left of $U$. Let $\delta' = L_U(\delta)$ be the unique isotropic Schur root in $\C(V', W')$. Since $M$ is Schur and not exceptional, $M$ has dimension vector $\delta'$.
\end{proof}

Let $E$ be an exceptional sequence of $\sigma$-stable representations. Clearly, such a sequence has length at most $n-2$. The sequence $E$ is said to be \emph{full} if $E$ has length $n-2$. The following result guarantees that such a full exceptional sequence of $\sigma$-stable representations always exists.

\begin{Lemma} \label{LemmaCompletion} Any exceptional sequence of $\sigma$-stable representations can be completed to a full exceptional sequence of $\sigma$-stable representations. In particular, there exists an exceptional sequence $(M_{n-2}, \ldots, M_2, M_{1})$ of $\sigma$-stable representations.
\end{Lemma}

\begin{proof} Let $(X_1, \ldots, X_r)$ be an exceptional sequence with all $X_i$ $\sigma$-stable. Assume that $r$ is not equal to $n-2$. Observe that $\langle d_{X_i}, d_{X_j} \rangle \le 0$ whenever $i \ne j$. By Proposition \ref{ExceptionalSequenceIsotropic}, and since the perpendicular category $\C(X_1, \ldots, X_r)^\perp$ contains representations of dimension vector $\delta$, we can extend it to an exceptional sequence
$$(X_1, \ldots, X_r, V,W)$$
of length less than $n$ where $\C(V,W)$ is tame and $\delta = d_V + d_W$. Therefore, there exists an exceptional representation $Y$ such that $(X_1, \ldots, X_r, Y, V, W)$ is exceptional. Being left orthogonal to both $V,W$, the representation $Y$ has to be $\sigma$-semistable. It is well known and easy to see that for any acyclic quiver, there exists a sincere representation of it that is rigid. Since the category $\C(X_1, \ldots, X_r,Y)$ is equivalent to the category of representations of an acyclic quiver with $r+1$ vertices, there exists a rigid object $Z \in \C(X_1, \ldots, X_r, Y)$ such that its Jordan-Hölder composition factors in $\C(X_1, \ldots, X_r, Y)$ will consist of all the simple objects in $\C(X_1, \ldots, X_r, Y)$. Being $\sigma$-stable, the objects $X_1, \ldots, X_r$ are non-isomorphic simple objects in $\C(X_1, \ldots, X_r, Y)$. Let $Y'$ be the other simple object. Being rigid, $Z$ is a general representation. It has some filtration where the subquotients are $\{X_1, \ldots, X_r, Y'\}$ with possible multiplicities. If $Y'$ is $\sigma$-stable, then we have a list $\{X_1, \ldots, X_r, Y'\}$ of $r+1$ objects that are $\sigma$-stable. Since these objects are the simple objects of $\C(X_1, \ldots, X_r, Y)$, they could be ordered to form an exceptional sequence of length $r+1$, which will be the wanted sequence. If $Y'$ is not $\sigma$-stable, then the above filtration of $Z$ can be refined to a filtration in $\rep(Q)$ where all subquotients are $\sigma$-stable. In particular, there will be more than $r$ non-isomorphic subquotients. Since $Z$ is rigid, we know that all subquotients have dimension vector a real Schur root, by Lemma \ref{Lemma0}. Now, these real Schur roots can be ordered to form an orthogonal sequence of real Schur roots; see \cite{DWSchurSeqn}. This sequence corresponds to an exceptional sequence of $\sigma$-stable representations of length at least $r+1$.
\end{proof}

\medskip

Now, we treat the case where $n$ is arbitrary.
We first need some more notations. From now on, let us fix  $(M_{n-2}, \ldots, M_{1})$ a full exceptional sequence of $\sigma$-stable representations.
Complete this sequence to get a full exceptional sequence $(M_{n-2}, \ldots, M_{1}, V, W)$. Note that $\delta$ is a root in $\C(V,W)$. We construct a sequence of isotropic Schur roots $\delta_1, \ldots, \delta_{n-2}$ and a sequence of rank three subcategories $\C_i$ as follows.  Set $\delta_1 = \delta$, $V_1 = V$ and $W_1=W$ and $\C_1 = \C(M_1, V_1, W_1)$. Observe that if $\C_1$ is wild, then it is connected and $M_1$ is preprojective or preinjective in $\C_1$. If $\C_1$ is tame, then it is either disconnected (in which case $M_1$ is the unique simple object in the trivial component of $\C_1$) or else it is connected and $M_1$ is quasi-simple in $\C_1$.
For $1 \le i \le n-2$, if $M_{i}$ is not preinjective in $\C_i=\C(M_{i}, V_i, W_i)$, then we reduce to the case of the exceptional sequence $(M_{n-2}, \ldots, M_{i+1}, V_{i+1}, W_{i+1})$ and we set $\delta_{i+1} = \delta_i, V_{i+1} = V_i, W_{i+1}=W_i$. If $M_{i}$ is preinjective in $\C(M_{i}, V_i, W_i)$, we reduce to the case of the exceptional sequence $(M_{n-2}, \ldots, M_{i+1}, V_{i+1}, W_{i+1})$ where $V_{i+1}$ is the reflection of $V_i$ to the left of $M_{i}$ and $W_{i+1}$ is the reflection of $W_i$ to the left of $M_{i}$. We set $\delta_{i+1} = \delta_i - \langle \delta_i, d_{M_{i}} \rangle d_{M_{i}}.$ In all cases, we set $\C_{i+1}=\C(M_{i+1}, V_{i+1}, W_{i+1})$.
Note that for $1 \le i \le n-2$, the root $\delta_i$ is a root in $\C(V_i, W_i)$. Finally, we set $\bar \delta = \delta_{n-2}$.

\medskip

For two dimension vectors $d_1, d_2$, we write $d_1 \hookrightarrow d_2$ provided a general representation of dimension vector $d_2$ has a subrepresentation of dimension vector $d_1$.

\begin{Lemma} We have $\delta_{i+1} \hookrightarrow \delta_{i}$ and each $\delta_{i}$ is an isotropic Schur root that is $\sigma$-semistable.
\end{Lemma}

\begin{proof} We proceed by induction on $i$. The case where $i=1$ is clear. Let $M(\delta_i)$ be a general representation of dimension vector $\delta_i$ in $\C(V_i, W_i)$ that is $\sigma$-semistable. If $M_{i}$ is not preinjective in $\C_i=\C(M_{i}, V_i, W_i)$, then $\delta_{i+1} = \delta_i$ and hence, it is clear that $\delta_{i+1} \hookrightarrow \delta_{i}$ and that $\delta_{i}$ is an isotropic Schur root that is $\sigma$-semistable. So assume that $M_{i}$ is preinjective in $\C_i$. Then $\Ext^1(M(\delta_i), M_{i})=0$.
If we have a non-zero morphism $f_i: M(\delta_i) \to M_{i}$ that is not an epimorphism, then the cokernel $C_i$ of $f_i$ is $\C_i$-preinjective and $\sigma$-semistable in $\C_i$. There is an epimorphism $C_i \to N_i$ where $N_i$ is relative $\sigma$-stable in $\C_i$. Since $M_i$ is $\C_i$-preinjective, $N_i$ (and $C_i$) has to be $\C_i$-preinjective and the only $\C_i$-preinjective relative $\sigma$-stable object in $\C_i$ is $M_{i}$. This is a contradiction. Therefore, any non-zero morphism $M(\delta_i) \to M_{i}$ is an epimorphism.
Now, take $e_{i} = \langle \delta_i, d_{M_{i}}\rangle = {\dim}\Hom(M(\delta_i), M_{i})$. We have a morphism $g_i: M(\delta_i) \to M_{i}^{e_{i}}$ given by a basis of $\Hom(M(\delta_i), M_{i})$, and as argued above,
$g_i$ is an epimorphism. The kernel $K_i$ of $g_i$ is of dimension $\delta_i - \langle \delta_i, d_{M_{i}}\rangle d_{M_{i}}$ which is $\delta_{i+1}$. Hence, $\delta_{i+1}$ is relative $\sigma$-semistable and hence $\sigma$-semistable. Consider the short exact sequence
$$0 \to K_i \to M(\delta_i) \to M_{i}^{e_{i}} \to 0$$
in $\C_i$. Since $K_i \in \hskip -2pt ^\perp M_{i}$ by construction, and $K_i$ has dimension vector $\delta_{i+1}$, we have $\ext(\delta_{i+1}, e_i\cdot d_{M_{i}})=0$ and by Schofield's result \cite[Theorem 3.3]{SchofieldGeneral}, $\delta_{i+1} \hookrightarrow \delta_{i+1} + e_i\cdot d_{M_{i}} = \delta_i$.
\end{proof}

Call $\delta$ of \emph{smaller type} if the $\tau$-orbit of a general representation of dimension vector $\delta$ contains a non-sincere representation. Equivalently, if $\tau$ denotes the Coxeter matrix of $Q$, then there is an integer $r$ such that $\tau^r \delta$ is not sincere.

\begin{Prop} \label{PropSmallerType}The following conditions are equivalent.
\begin{enumerate}[$(1)$]
    \item The root $\delta$ is of smaller type.
\item There is a $\sigma$-stable representation that is preprojective or preinjective.
\item There is a $\sigma$-semistable representation that is preprojective or preinjective.\end{enumerate}
\end{Prop}

\begin{proof}
Let $Z$ be a $\sigma$-semistable representation that is preprojective or preinjective. If $Z$ is preprojective, then there is some $i \ge 0$ with $\tau^iZ=P$ projective. Since $\langle d_Z, \delta \rangle=0$, we get $\langle d_P, \tau^i\delta \rangle=0$, showing that $\tau^i\delta$ is not sincere, that is, $\delta$ is of smaller type. If $Z$ is preinjective, then there is some $j < 0$ with $\tau^jZ=Q[1]$ a shift of a projective representation $Q$. Since $\langle d_Z, \delta \rangle=0$, we get $\langle -d_Q, \tau^j\delta \rangle=0$, showing that $\tau^j\delta$ is not sincere, that is, $\delta$ is of smaller type. This proves that $(3)$ implies $(1)$. Clearly, $(2)$ implies $(3)$. Assume that $\delta$ is of smaller type. Then there is some integer $i$ and some $x \in Q_0$ such that $\Hom(P_x, \tau^iM(\delta))=0$ whenever $M(\delta)$ has dimension vector $\delta$. Therefore, there exists a preprojective (if $i \ge 0$) or preinjective (if $i < 0$) representation $Z$ such that $Z$ is left orthogonal to any representation $M(\delta)$ of dimension vector $\delta$. Therefore $Z$ is $\sigma$-semistable. If $Z$ is not $\sigma$-stable and preprojective, then it has a $\sigma$-stable subrepresentation $Z'$ which has to be preprojective. If $Z$ is not $\sigma$-stable and preinjective, then it has a $\sigma$-stable quotient $Z'$ which has to be preinjective.
\end{proof}

We start with the following result, which describes the $\sigma$-stable objects when $n=4$. The core of it will be generalized to arbitrary $n$ later. We think it is interesting to have a separate result for the $n=4$ case since more can be said for this small case.

\begin{Prop} \label{Prop4Vertices}
Let $Q$ be a connected wild quiver with $4$ vertices and let $(M_2, M_1)$ be a full exceptional sequence of $\sigma$-stable representations. Consider the exceptional sequence $(M_2, M_1, V, W)$ with $V,W$ the simple objects in $\C(V,W)$.  Let $M$ be $\sigma$-stable but not isomorphic to $M_1$ or $M_2$. Then
\begin{enumerate}[$(1)$]
\item The root $\bar \delta$ is the only $\sigma$-stable non-real Schur root.
\item The root $\delta$ is of smaller type if and only if one of $M_1, M_2$ is preprojective or preinjective.
\item If $\delta$ is not of smaller type, then the only exceptional $\sigma$-stable representations are $M_1, M_2$.
\end{enumerate}
\end{Prop}

\begin{proof}
The sufficiency of (2) follows from Proposition \ref{PropSmallerType}.
Assume now that both $M_1, M_2$ are regular. We claim that $\delta$ is not of smaller type and that the only exceptional $\sigma$-stable representations are $M_1, M_2$. In this case, by Lemma \ref{ResultBaerStrauss}, since $\C_1 =\C(M_1, V, W)$ is wild and $\C(V, W)$ is tame, we know that $M_1$ is either preprojective or preinjective in $\C_1$. Observe that $\C_1$ is connected since it is wild and $\delta$ is an isotropic root in $\C_1$. If $M_1$ is preinjective in $\C_1$, then $\C_2 = ^\perp M_1$ is also wild, since $M_1$ is regular. If $M_1$ is preprojective in $\C_1$, consider the exceptional sequence $(M_1', M_2, V, W) = (M_1', M_2, V_2, W_2)$. Note that $M_1, M_2$ are the relative simples in $\C(M_2,M_1)$. Since any indecomposable object in $\C(M_2,M_1)$ has a morphism from and to $M_1 \oplus M_2$, all objects of $\C(M_2,M_1)$, seen as objects in $\rep(Q)$, are regular. Therefore, $M_1' \in \C(M_2,M_1)$ is regular in $\rep(Q)$ and $M_1'^\perp = \C(M_2, V_2, W_2)=\C_2$ is wild (and connected). Let $M$ be an arbitrary $\sigma$-stable representation not isomorphic to $M_1$ or $M_2$. Since $\Hom(M_2, M)=0$, there exists a non-negative integer $d_2$ and a short exact sequence
$$(*): \quad 0 \to M \to E_2 \stackrel{f_1}{\to} M_2^{d_2} \to 0$$
with $E_2 \in M_2^\perp$. Observe that $d_2 = 0$ if and only if $M \in M_2^\perp$. Observe also that $\Hom(E_2, M_1) = \Hom(M_1, E_2)=0$.
Since $\Hom(M_1, E_2)=0$, there exists a non-negative integer $d_1$ and a short exact sequence
$$(*): \quad 0 \to E_2 \to E_1 \stackrel{f_1}{\to} M_1^{d_1} \to 0$$
with $E_1 \in M_2^\perp \cap M_1^\perp$. Observe that $d_1 = 0$ if and only if $E_2 \in M_1^\perp$. Since $E_1$ is relative $\sigma$-semistable in $M_2^\perp \cap M_1^\perp=\C(V,W)$, there is a monomorphism $M(\delta) \to E_1$ where $M(\delta)$ is a Schur representation of dimension vector $\delta$ in $\C(V,W)$.  If $M_1$ is preinjective in $\C_1$, then $\Ext^1(M(\delta), M_1)=0$. Then $e_1:=\langle \delta, d_{M_1} \rangle = {\rm dim}\Hom(M(\delta), M_1)>0$. We have a short exact sequence
$$0 \to M(\delta_1) \to M(\delta)\to M_1^{e_1}\to 0$$
and this yields a monomorphism $M(\delta_1) \to E_2$. If $M_1$ is preprojective in $\C_1$, then $\Hom(M(\delta), M_1)=0$ and hence, we get a monomorphism $M(\delta_1) = M(\delta) \to E_2$. If $M_2$ is preinjective in $\C_2$, then $\Ext^1(M(\delta_1), M_2)=0$. Then $e_2:=\langle \delta_1, d_{M_2} \rangle = {\rm dim}\Hom(M(\delta_1), M_2)>0$. We have a short exact sequence
$$0\to M(\delta_2) \to M(\delta_1)\to M_2^{e_2}\to 0$$
and this yields a monomorphism $M(\delta_2) \to M$. If $M_2$ is preprojective in $\C_2$, then $\Hom(M(\delta_1), M_2)=0$ and hence, we get a monomorphism $M(\delta_2) = M(\delta_1) \to M$. This shows that $M \cong M(\delta_2)=M(\bar\delta)$. In particular, the exceptional sequence $(M_2,M_1)$ of $\sigma$-stable representations is unique. It follows from Proposition \ref{PropSmallerType} that $\delta$ is not of smaller type. This prove our claim and hence $(2)$ and $(3)$.
Statement $(1)$ in the case where $\delta$ is of smaller type is a consequence of the next proposition, Proposition \ref{MainProp}.
\end{proof}

Given a subcategory $\C$ of $\rep(Q)$, we denote by $\tau_\C$ the Auslander-Reiten translation in $\C$. Fix a full exceptional sequence $(M_{n-2}, \ldots, M_2, M_1)$ of $\sigma$-stable representations, and complete it to form an exceptional sequence $(M_{n-2}, \ldots, M_2, M_1, V, W)$ where $V,W$ are the non-isomorphic simple objects in $\C(V,W)$. In particular, $\delta = d_V + d_W$ and $\langle d_W, d_V \rangle = -2$. Recall that for $1 \le i \le n-2$, we denote by $\C_i$ the subcategory $\C(M_i, V_i, W_i)$, where $V_i, W_i$ have been defined previously. If $\C_i$ is tame and connected, then $M_i$ is regular quasi-simple in $\C_i$ and lies in a tube of rank $2$ in $\C_i$. We set $N_i := \tau_{\C_i}M_i = \tau_{\C_i}^{-1}M_i$. The exceptional sequence $(M_{n-2}, \ldots, M_{i+1})$ is denoted by $\mathcal{E}_i$.

\begin{Prop} \label{MainProp}
Let $Q$ be a connected quiver with isotropic Schur root $\delta$. Let $N$ be a $\sigma$-stable object not isomorphic to any $M_i$.
\begin{enumerate}[$(1)$]
    \item The root $\bar \delta$ is the only $\sigma$-stable Schur root that is not real.
\item If $N$ is not exceptional, then $d_N = \bar \delta$ and $N \in \C(V_{n-2},W_{n-2})$.
\item If $N$ is exceptional, then there exists some $i$ such that $\C_i$ is tame connected with $M_i$ quasi-simple. There exists a subsequence $\mathcal{F}_i$ of $\mathcal{E}_i$ such that $N$ is the reflection of $N_i$ to the left of $\mathcal{F}_i$.\end{enumerate}
\end{Prop}

\begin{proof}
We first need to check that there is at least one $\sigma$-stable root that is not real. Assume otherwise. As $\delta$ is $\sigma$-stable, its $\sigma$-stable decomposition will involve only real Schur roots. These roots may be ordered to form a Schur sequence. Since all roots are real, this will correspond to an exceptional sequence $(F_1, \ldots, F_r)$. Suppose that $r < n$. Complete the latter sequence to get a full exceptional sequence $(F_1, \ldots, F_n)$. Since $\delta$ is a root in each $F_i^\perp$ for $1 \le i \le r$, we see that $\delta$ is a root in $\C(F_{r+1}, \ldots, F_n)$. Since $\delta$ is also a root in $\C(F_{1}, \ldots, F_r)$, we get that $\delta$ is in the span of $d_{F_1}, \ldots, d_{F_r}$ as well as in the span of  $d_{F_{r+1}}, \ldots, d_{F_n}$. Since $(F_1, \ldots, F_n)$ is an exceptional sequence, the vectors  $d_{F_1}, \ldots, d_{F_n}$ are linearly independent, a contradiction. Hence, $r=n$. But then, $\C(F_{1}, \ldots, F_n) = \rep(Q)$ and all objects of $\rep(Q)$ are $\sigma$-semistable, which is also a contradiction.
In the rest of the proof, we will see that the only possibility for a non-real $\sigma$-stable Schur root is $\bar \delta = \delta_{n-2}$.

Let $N$ be a $\sigma$-stable representation not isomorphic to any $M_i$. In particular, $\Hom(N, M_i) = \Hom(M_i,N)=0$ for all $i$. Let $E_{n-2}:=N$. For $i=n-2, n-3, \ldots, 2, 1$, we have a short exact sequence
$$0 \to E_i \to E_{i-1} \to M_{i}^{d_{i}}\to 0$$
where $d_{i}$ is a non-negative integer with $d_{i} = -{\rm dim}\Ext^1(M_{i},E_i) = \langle d_{M_{i}}, d_{E_i} \rangle \ge 0$, as $\Hom(M_i, E_i)=0$ by induction. Observe that for $i < n-2$, we have $E_i \in M_{n-2}^\perp \cap \cdots \cap M_{i+1}^\perp$. In particular, $E_0$ lies in $\C(V,W)$. Observe also that all $E_i$ are $\sigma$-semistable.

Since $E_0$ is relative $\sigma$-semistable in $\C(V,W)$, there is a monomorphism $Z_0:=M(\delta) \to E_0$ where $M(\delta)$ is a Schur representation of dimension vector $\delta$ in $\C(V,W)$.  If $M_1$ is preinjective in $\C_1$, then $\Ext^1(M(\delta), M_1)=0$. Then $e_1:=\langle \delta, d_{M_1} \rangle = {\rm dim}\Hom(M(\delta), M_1)$. We have a short exact sequence
$$0 \to M(\delta_2) \to M(\delta)\to M_1^{e_1}\to 0$$
and this yields a monomorphism $M(\delta_2) \to E_1$ where $M(\delta_2)$ lies in $\C(V_2, W_2)$ is $\sigma$-semistable. If $M_1$ is preprojective in $\C_1$, then $\Hom(M(\delta), M_1)=0$ and hence, we get a monomorphism $M(\delta_2) = M(\delta) \to E_1$. If $\C_1$ is tame disconnected, then $\Hom(M(\delta), M_1)=0$ and we get a monomorphism $M(\delta_2) = M(\delta) \to E_1$ as well. Assume that $\C_1$ is tame connected. If $\Hom(M(\delta), M_1)=0$, then we get a monomorphism $M(\delta_2) = M(\delta) \to E_1$ as previously. Assume that $\Hom(M(\delta), M_1)\ne 0$. Then $\Hom(N_1, M_1)=0$ and there is a monomorphism $N_1 \to M(\delta)$. Therefore, we get a monomorphism $N_1 \to E_1$. Hence, in all cases, we have a monomorphism $Z_1 \to E_1$ where $Z_1$ is either a Schur ($\sigma$-semistable) representation of dimension vector $\delta_2$ in $\C(V_2, W_2)$ or is $N_1$. In the first case, we consider the exceptional sequence $(M_{n-2}, \ldots, M_2, V_2, W_2)$ and we proceed by induction starting with the monomorphism $Z_1 \to E_1$ where $Z_1 = M(\delta_2)$. For the second case, we consider the exceptional sequence $(M_{n-2}, \ldots, M_2, N_1)$. If $\Hom(N_1, M_2) \ne 0$, then $\Ext^1(N_1, M_2) = 0$ as $\Ext^1(M_2, N_1) = 0$. We have a short exact sequence
$$0 \to Z_2 \to N_1 \to M_2^{f_2} \to 0$$
where $f_2 = {\rm dim} \Hom(N_1, M_2)$. Then, $Z_2 \in ^\perp M_2$ and hence, we get a monomorphism $Z_2 \to E_2$ and $Z_2$ is the reflection of $N_1$ to the left of $M_2$. If $\Hom(N_1, M_2) = 0$, then we get a monomorphism $Z_2:=N_1 \to E_2$. We proceed by induction until we get a monomorphism $Z_{n-2} \to E_{n-2}=N$. This gives $Z_{n-2} \cong N$. We see that $Z_{n-2}$ will be of the required form. In particular, if $N$ is not exceptional, then $N$ has dimension vector $\delta_{n-2}=\bar \delta$ and will be in $\C(V_{n-2}, W_{n-2})$.
\end{proof}

The next two lemmas will help us in giving a better description of the simple objects in $\A(\delta)$, that is, the $\sigma$-stable objects.

\begin{Lemma} \label{Lemma4}
Assume that both $\C(M_{t-1}, \ldots, M_1, V, W)$ and $\C_t$ are tame connected with $\delta$ the unique isotropic Schur root in both of these subcategories. Then $\C(M_t, M_{t-1}, \ldots, M_1, V, W)$ is tame connected with only one isotropic Schur root $\delta$.
\end{Lemma}

\begin{proof}Assume that $\C:=\C(M_t, M_{t-1}, \ldots, M_1, V, W)$ is wild. Then $M_t$ lies in the preprojective or preinjective component of a wild connected component of $\C$. Observe that, in $\C_t$, the object $M_t$ is left orthogonal to the unique isotropic Schur root (which is $d_{V_t} + d_{W_t}$). Since $\C_t$ is connected, this means that $M_t$ is regular in $\C_t$. In particular, there are infinitely many indecomposable objects $Z$ of $\C_t$ (and hence of $\C$) with $\Hom(M_t, Z)\ne 0$. Therefore, $M_t$ cannot lie in a preinjective component of $\C$. Similarly, $M_t$ cannot lie in a preprojective component of $\C$. This is a contradiction. It remains to show that $\C$ is connected (it is clear that then, $\delta$ will be the unique isotropic Schur root in $\C$). Assume otherwise. Then, we have $\C \cong \mathcal{B}_1 \times \mathcal{B}_2$ where each of $\mathcal{B}_1, \mathcal{B}_2$ is equivalent to a category of representations of a non-empty acyclic quiver. Assume that $\delta \in \mathcal{B}_2$. By assumption, we must have that both $\C(M_{t-1}, \ldots, M_1, V, W), \C_t$ are subcategories of $\mathcal{B}_2$. But then, $\C = \mathcal{B}_2$, a contradiction. Therefore, $\C$ is connected.
\end{proof}

\begin{Lemma} \label{Lemma5} Assume that $\C_i$ is tame connected for all $1 \le i \le s-1$. Fix $1 \le t \le s-1$ and consider the representation $N_t$ as defined above.
\begin{enumerate}[$(1)$]
    \item If $M_s$ is preprojective in the category $\C_s$, then $\Hom(N_t,M_s)=0$.
\item If $M_s$ is preinjective in the category $\C_s$, then $M_i \in ^\perp M_s$ for all $1 \le i \le s-1$.
\item If $M_s$ is simple disconnected in $\C_s$, then all $\sigma$-stable representations, except possibly $M_{s-1}, \ldots, M_1$, lie in $^\perp M_s$.
\end{enumerate}
\end{Lemma}

\begin{proof} Consider the categories $$\C := \C(M_s, \ldots, M_1, V, W)$$ and $$\C':=\C(M_{s-1}, \ldots, M_1, V, W).$$ We know that $\C'$ is tame connected by Lemma \ref{Lemma4}.
Assume that $\C$ is wild. Therefore, $M_s$ lies in a preprojective or preinjective component of $\C$. For proving (1), assume that $M_s$ is preprojective in the category $\C_s$, which means that $M_s$ is preprojective in $\C$. Assume to the contrary that $\Hom(N_t,M_s)\ne 0$. Since $N_t \in \C$, the object $N_t$ is preprojective in $\C$. Therefore, $N_t$ cannot be regular in $\C_t$, a contradiction. For proving (2), assume that $M_s$ is preinjective in the category $\C_s$, which means that $M_s$ is preinjective in $\C$. Since $t$ is arbitrary and $\Hom(M_t, M_i)=0$, it is sufficient to prove that $\Ext^1(M_t, M_s)=0$. If not, then the Auslander-Reiten formula in $\C$ yields a nonzero morphism from $M_s$ to $\tau_\C M_t$ and hence, $\tau_\C M_t$ (and thus $M_t$) is preinjective in $\C$. Thus, $M_t$ cannot be regular in $\C_t$, a contradiction. Clearly, if $\C$ is tame, then (1), (2) cannot occur, since a subcategory of a tame category cannot be wild.

It remains to prove (3). Let $N$ be a $\sigma$-stable object not isomorphic to any $M_i$. It follows from the proof of the last theorem that we have a short exact sequence
$$0 \to E_s \to E_{s-1} \to M_s^{d_s}\to 0$$ and a monomorphism $Z \to E_{s-1}$ where $Z$ is either indecomposable of dimension vector $\delta$ in $\C(V,W)$ or is exceptional and relative $\sigma$-semistable in $\C'$ (a reflection of one of the $N_i$ for $1 \le i \le s-1$). In the first case, since $\C_s$ is tame disconnected, this yields a monomorphism $Z \to E_s$ and the rest of the proof of the above theorem deals with representations in $^\perp M_s$. In the second case, we get a monomorphism $M \to Z$ where $M$ is relative $\sigma$-stable in $\C'$ and exceptional, and hence a quasi-simple object of $\C'$. Assume that $M$ lies in a tube $T$ of rank $r$. In particular, the other non-isomorphic $r-1$ quasi-simple objects of that tube are among the objects $M_{s-1}, \ldots, M_1$. We claim that $\Ext^1(M,M_s)=0$. There exists a short exact sequence
$$0 \to M \to Y \to N \to 0$$ in $\C'$ where $M,Y,N$ are indecomposable in $T$ with $Y$ of dimension vector $\delta$ and $N$ of quasi-length $r-1$ which does not have $M$ as a quasi-simple composition factor. Using that $(M_{s-1}, \ldots, M_1, V,W)$ is an exceptional sequence and since, by construction, $M$ is the unique quasi-simple of $T$ with $\Hom(M,Y)\ne 0$, we have  $Y \in \C(V,W)$. The surjective map $\Ext^1(Y,M_s) \to \Ext^1(M,M_s)$ together with the fact that $\C_s = \C(M_s, V, W)$ is disconnected gives that $\Ext^1(M,M_s)=0$. This proves our claim.  Observe that $M$ lies in $M_s^\perp$. We have a short exact sequence
$$0 \to Z' \to M \to M_s^{f} \to 0$$
where $f = {\rm dim} \Hom(M, M_s) = \langle d_M, d_{M_s} \rangle$. Thus, $Z'$ is the reflection of $M$ to the left of $M_s$ and $Z' \in ^\perp M_s$ and we get a monomorphism $Z' \to E_s$. The proof of the last theorem continues with $Z_s = Z'$ and the monomorphism $Z_s \to E_s$. In all cases, we see that $M_s$ satisfies the required property.
\end{proof}

\section{Cone of $\sigma$-semi-stable dimension vectors}

Let $d$ be a dimension vector. Let us denote by $C(\sigma_d)$ the set of all $\sigma_d$-semistable dimension vectors. We consider $C_\mathbb{R}(\sigma_d)$ the corresponding cone in $\mathbb{R}^n$, which lie in the positive orthant of $\mathbb{R}^n$. Since the rays $x \in C_\mathbb{R}(\sigma_d)$ satisfy $\langle x, d \rangle = 0$, we rather consider $C_\mathbb{R}(\sigma_d)$ as a cone in $\mathbb{R}^{n-1}$. The integral vectors in $C_\mathbb{R}(\sigma_d)$ correspond to the $\sigma_d$-semistable dimension vectors. For $d=\delta$, since there exists a full exceptional sequence $(M_{n-2}, \ldots, M_2, M_1)$ of $\sigma_\delta$-stable representations and since the dimension vectors in an exceptional sequence are linearly independent, we see that $C_\mathbb{R}(\sigma_\delta)$ is a cone of full dimension in $\mathbb{R}^{n-1}$. In general, it is well known that $C_\mathbb{R}(\sigma_d)$ is a cone over a polyhedron where the indivisible dimension vectors in the extremal rays are $\sigma_d$-stable dimension vectors. On the other hand, a $\sigma_d$-stable dimension vector needs not lie on an extremal ray.

\begin{Lemma}
If $f$ is $\sigma_d$-stable and a real Schur root, then it lies in an extremal ray.
\end{Lemma}

\begin{proof}
Let $f$ be a $\sigma_d$-stable real Schur root. Assume that $f$ does not lie on an extremal ray. According to \cite[Theorem 6.4]{DWSchurSeqn}, there exists dimension vectors $f_1, \ldots, f_s$, all lying on extremal rays, such that $f$ is a positive integral combination of $f_1, \ldots, f_s$ and $f_1, \ldots, f_s$ are linearly independent in $\mathbb{R}^n$. Moreover, we have $\langle f, f_i \rangle \le 0$ and $\langle f_i, f \rangle \le 0$ for all $1 \le i \le s$. These conditions imply that $$1=\langle f, f \rangle = \langle f, a_1f_1 + \cdots + a_sf_s \rangle \le 0,$$ a contradiction.
\end{proof}

There is a special case of interest, which is when there is some ray $[r]$ such that all $\sigma_d$-stable dimension vectors, except possibly the ones on $[s]$, lie on extremal rays (this is the case when $d=\delta$ where $\delta$ is an isotropic Schur root). In such a case, either all $\sigma_d$-stable dimension vectors lie on the boundary of the cone $C_\mathbb{R}(\sigma_d)$ or else $[s]$ lies in the interior of $C_\mathbb{R}(\sigma_d)$.

\medskip

Let $v_1, v_2, \ldots, v_r$ be the extremal rays of $C_\mathbb{R}(\sigma_d)$. Take $f$ any dimension vector lying in $C_\mathbb{R}(\sigma_d)$ but that is neither in an extremal ray nor in the ray $[s]$. Then we know that $f$ has at least one $\sigma_d$-stable factor that lies on an extremal ray.
Since the $\sigma_d$-stable factors of $f$ can be ordered to form an orthogonal sequence of Schur roots, we see that there exists a Schur root $\alpha$ (that corresponds to an extremal ray $v_i$) such that either $\langle v_i, f \rangle > 0$ or $\langle f, v_i \rangle > 0$. If $\langle v_i, f \rangle > 0$, then $\langle v_i, x \rangle = 0$ defines an hyperplane cutting $C_\mathbb{R}(\sigma_d)$ such that $f, v_i$ lie on the same side while all other extremal rays $v_1, \ldots, v_{i-1}, v_{i+1}, \ldots, v_r$ of $C_\mathbb{R}(\sigma_d)$ lie on the other side or on the boundary of that hyperplane. We get a similar situation if $\langle f, v_i \rangle > 0$ by considering the hyperplane $\langle x, v_i \rangle = 0$.

\medskip

For $1 \le i \le r$, let $C_i(\sigma_d)$ be the cone in $\mathbb{R}^{n-1}$ generated by all the rays $v_1, \ldots, v_r$ but $v_i$. By the above observation, we have that $${\rm Proper}(C_\mathbb{R}(\sigma_d)):=\cap_{1 \le i \le r}C_i(\sigma_d)$$ is either empty or else contains only the ray $[s]$. We will see that this restriction yield a very beautiful description of $C_\mathbb{R}(\sigma_d)$ and, in particular, of $C_\mathbb{R}(\sigma_\delta)$.

\begin{Prop} \label{PropConeSimplex}
We have that ${\rm Proper}(C_\mathbb{R}(\sigma_d)) = \emptyset$ if and only if $C_\mathbb{R}(\sigma_d)$ is the cone over a simplex.
\end{Prop}

\begin{proof} The sufficiency is easy to see. Assume that ${\rm Proper}(C_\mathbb{R}(\sigma_d)) = \emptyset$. We may work in $\mathbb{R}^{j}$ where $j \le n-1$ and assume that $C_\mathbb{R}(\sigma_d)$ is of full dimension in $\mathbb{R}^{j}$. By Radon's theorem, if the number of extremal rays $r$ of $C_\mathbb{R}(\sigma_d)$ is at least $(j-1)+2 = j+1$, then we can partition the rays $v_1, \ldots, v_r$ into two non-empty subsets $A,B$ such that the corresponding cones $C(A)$ and $C(B)$ generated by the rays in $A$ and by the rays in $B$ have a ray of intersection. This ray of intersection will have to be in ${\rm Proper}(C_\mathbb{R}(\sigma_d)) = \emptyset$, a contradiction. Therefore, $r \le j$. Since $C_\mathbb{R}(\sigma_d)$ is of full dimension in $\mathbb{R}^{j}$, then $r=j$ and $C_\mathbb{R}(\sigma_d)$ is the cone over an $(j-1)$-simplex.
\end{proof}

Now, we are interested in the case where ${\rm Proper}(C_\mathbb{R}(\sigma_d))$ is reduced to a single ray $[s]$ (which then has to be the ray of $\bar \delta$ if $d = \delta$). Let us take an affine slice $\Delta$ of $C_\mathbb{R}(\sigma_d)$. The rays $v_1, \ldots, v_r$ will correspond to points $u_1, \ldots, u_r$ in $\Delta$ and these points are the vertices of a polyhedron $\mathcal{R}$ in $\Delta$ defined as the convex hull of $u_1, \ldots, u_r$. The ray $[s]$ corresponds to a point $s$ in $\mathcal{R}$. In order to study the convex properties of $\mathcal{R}$, let us translate $\mathcal{R}$ so that $s$ coincides with the origin. In other words, set $w_i = u_i - s$ and consider the polyhedron $\mathcal{P}$ which is the convex hull of $w_1, \ldots, w_r$. Since $\mathcal{R}$ lies on an affine slice, we see that $\mathcal{P}$ lies in a subspace of dimension $n-2$ of $\mathbb{R}^{n-1}$.
Let $\mathcal{P}_i$ be the convex hull of all points $w_1, \ldots, w_r$ but $w_i$. We define $${\rm Proper}(\mathcal{P}):=\cap_{1 \le i \le r}\mathcal{P}_i$$
and we will be interested in the case where ${\rm Proper}(\mathcal{P})$ only contains the origin.

\medskip

The first two lemmas are easy to prove.

\begin{Lemma} \label{Lemma2}Let $\mathcal{P}'$ be the convex hull of a subset of $w_1, \ldots, w_r$. Then ${\rm Proper}(\mathcal{P}')$ is either empty or reduced to the origin.
\end{Lemma}

\begin{Lemma} \label{Lemma3}
Consider a non-trivial partition $\{w_{i_1}, \ldots, w_{i_s}\} = A_1 \cup A_2$ of a subset of $\{w_1, \ldots, w_r\}$. Denote by $\mathcal{P}_{A_i}$ the convex hull of the points in $A_i$, for $i=1,2$. Then  $\mathcal{P}_{A_1} \cap \mathcal{P}_{A_2}$ is either empty or reduced to the origin.
\end{Lemma}

\begin{Prop} \label{PropVectSpacesDecomp} Suppose that ${\rm Proper}(\mathcal{P})$ is empty or reduced to the origin and is full dimensional in $\mathbb{R}^{t}$ with $t \le n-2$. Then there exists a vector space decomposition $$\mathbb{R}^{t} = V_1 \oplus \cdots \oplus V_s$$ of $\mathbb{R}^{t}$ such that if $V_i$ has dimension $d_i$, then it contains $d_i + 1$ points among $0, w_1, \ldots, w_r$ that form a $d_i$-simplex in $V_i$ containing the origin.
\end{Prop}

\begin{proof}
We may assume that $t = n-2$ so that $\mathcal{P}$ is $(n-2)$-dimensional. If ${\rm Proper}(\mathcal{P})$ is empty, then $s=1$ and the result follows from Proposition \ref{PropConeSimplex}. Assume that ${\rm Proper}(\mathcal{P})$ is reduced to the origin. Suppose first that the origin lies on a facet, say $F$, of $\mathcal{P}$. We claim that $F$ contains $r-1$ of the points $w_1, \ldots, w_r$. Assume otherwise. Consider an $(n-3)$-simplex in $F$ generated by points $w_{i_1}, \ldots, w_{i_{n-2}}$. Let $u, v \in \{w_1, \ldots, w_r\}$ be two distinct points not in $F$. By Radon's theorem, we can partition the points $\{w_{i_1}, \ldots, w_{i_{n-2}}, u, v\}$ into two non-empty subsets $A_1, A_2$ such that $\mathcal{P}_{A_1} \cap \mathcal{P}_{A_2} \ne \emptyset$, where $\mathcal{P}_{A_i}$ denotes the convex hull of the points in $A_i$. By Lemma \ref{Lemma3}, this intersection is the origin and hence lies on $F$. Since $F$ is a facet, $u,v$ lie on the same side of $F$. Therefore, for $i=1,2$, the set $B_i:=A_i\backslash\{u,v\}$ is not empty. Now, $B_1$ and $B_2$ form a partition of $w_{i_1}, \ldots, w_{i_{n-2}}$ such that $\mathcal{P}_{B_1} \cap \mathcal{P}_{B_2} \ne \emptyset$, where $\mathcal{P}_{B_i}$ denotes the convex hull of the points in $B_i$. This contradicts Proposition \ref{PropConeSimplex}.

\medskip

Now, let us assume that the origin lies in the interior of $\mathcal{P}$. By Radon's theorem, we can write $\{w_1, \ldots, w_r\} = E_1 \cup E_2$ where $E_1,E_2$ are disjoint and non-empty such that $\mathcal{P}_{E_1} \cap \mathcal{P}_{E_2} = \{0\}$, where $\mathcal{P}_{E_i}$ denotes the convex hull of the points in $E_i$. By Carathéodory's theorem, there is a simplex formed by some points $z_1, \ldots, z_s$ in $E_1$ that contains the origin in its interior. With no loss of generality, assume that $z_i = w_i$ and $s \le r-1$. Let $V_1$ be the vector space spanned by the points $w_1, \ldots, w_s$ and consider the vector space $V_2$ spanned by the points $w_{s+1}, \ldots, w_r$. Let $C_1$ be the convex hull of the points $w_1, \ldots, w_s$ and let $C_2$ be the convex hull of the points $w_{s+1}, \ldots, w_r$. Since both $\mathcal{P}_{E_1}, \mathcal{P}_{E_2}$ contain the origin, we see that $C_1 \cap C_2 = \{0\}$.
Since $0$ lies in the interior of $C_1$, we get also that $V_1 \cap C_2=0$. We claim that $V_1 \cap V_2 = 0$. Assume that $V_1 \cap V_2$ is non-zero. Observe that any element in $V_1$ can be written as a non-negative linear combination of $w_1, \ldots, w_s$. There exists non-negative real numbers $a_1, \ldots, a_s$ and real numbers $b_{s+1}, \ldots, b_r$ such that
$$a_1w_1 + \cdots + a_sw_s = b_{s+1}w_{s+1} +  \cdots + b_rw_r.$$
Moreover, $a_1w_1 + \cdots + a_sw_s$ is non-zero. We may assume the $a_i$ small enough so that the left-hand side lies in $C_1$. Let us write $\{s+1, \ldots, r\} = I_1 \cup I_2$ where $I_1, I_2$ are disjoint and $i \in I_1$ if and only if $b_i \ge 0$. We may assume further that the $|b_i|$ are small enough so that both $\sum_{i \in I_1}b_iw_i, -\sum_{j \in I_2}b_jw_j$ lie in $C_2$. If $I_2 = \emptyset$, then $a_1w_1 + \cdots + a_sw_s \in C_1 \cap C_2$ is non-zero, a contradiction. If all $b_i$ are non-positive, then
$-\sum_{j \in I_2}b_jw_j$ lies in $V_1\cap C_2=\{0\}$, a contradiction. If some $b_i$ are negative and some $b_i$ are positive, we can rewrite the sum as
$$a_1w_1 + \cdots + a_sw_s + -\sum_{j \in I_2}b_jw_j= \sum_{j \in I_1}b_jw_j.$$
Considering Lemma \ref{Lemma3} with the partition $$(\{w_1, \ldots, w_s\}\cup\{w_{i}\mid i \in I_2\}) \cup (\{w_{j} \mid j \in I_1\})$$ of $\{w_1, \ldots, w_r\}$, we get that $a_1w_1 + \cdots + a_sw_s +  -\sum_{j \in I_2}b_jw_j$ is zero, which reduces to a case we have already considered. Therefore, we have proven that $\mathbb{R}^{n-2} = V_1 \oplus V_2$, where $V_1$ satisfies the property of the statement. We proceed by induction on $V_2$ with the points $w_{s+1}, \ldots, w_r$ and by using Lemma \ref{Lemma2}.
\end{proof}

\section{The ring of semi-invariants of an isotropic Schur root}

In this section, we denote by $\delta$ an isotropic Schur root and by $\sigma = \sigma_\delta$ the weight given by $-\langle -, \delta \rangle$.

Consider, as previously, a full exceptional sequence $(M_{n-2}, \ldots, M_2, M_1, V,W)$ where $(M_{n-2}, \ldots, M_2, M_1)$ is an exceptional sequence of simple objects in $\A(\delta)$. Take $I \subseteq \{1, \ldots, n-2\}$ such that $i \in I$ if and only if $\C_i$ is tame connected.

\begin{Defn}  The \emph{associated tame subcategory of $Q$ relative to $\delta$}, denoted $\mathcal{R}(Q,\delta)$, is the thick subcategory of $\rep(Q)$ generated by $(\bigoplus_{i \in I}M_i)\oplus V_{n-2}\oplus W_{n-2}$.
\end{Defn}

\begin{Theo} \label{MainTheo}Let $Q$ be an acyclic connected quiver and $\delta$ an isotropic Schur root. Then
\begin{enumerate}[$(1)$]
    \item The category $\mathcal{R}(Q,\delta)$ is tame connected with isotropic Schur root $\bar \delta$ and is uniquely determined by $\delta$.
    \item The simple objects in $\A(\delta)$, up to isomorphism, are given by the disjoint union $$\{M_i \mid i \not \in I\} \cup \{\text{quasi-simple objects in} \; \mathcal{R}(Q,\delta)\}.$$
    \item We have $${\rm SI}(Q, \delta) \cong {\rm SI}(\mathcal{R}(Q,\delta), \bar\delta)[x_{r+1}, \ldots, x_n].$$
\end{enumerate}
\end{Theo}

\begin{proof}Let $(M_{n-2}, \ldots, M_1)$ be an exceptional sequence of $\sigma_\delta$-stable representations with the corresponding full exceptional sequence $(M_{n-2}, \ldots, M_1, V, W)$ in $\rep(Q)$, where $\delta = d_V + d_W$. First, denote by $M_{l_r}, \ldots, M_{l_1}$ with $l_r > \cdots > l_1$ the $M_j$ such that $\C_j$ is tame disconnected or such that $M_j$ is preprojective in $\C_j$. We get an exceptional sequence $$(*) \qquad (N_{t}, \ldots, N_2, N_1, V,W, M_{l_r}', \ldots, M_{l_1}')$$ where all $M_{l_j}$ have been reflected, one by one, to the right of the exceptional sequence. Observe that $\{M_i \mid i \in I\} \subseteq \{N_1, \ldots, N_t\}$. Let
$$ \{N_t, \ldots, N_1\} \backslash \{M_i \mid i \in I\}:=\{N_{j_s},\ldots, N_{j_1}\}$$
where $j_s > \cdots > j_1$. Assume also that $I = \{m_1, \ldots, m_q\}$ with $m_q > \cdots > m_1$. We have $q+s=t$ and $t+r=n-2$. By Lemma \ref{Lemma5}(2), we may reflect all exceptional objects of $\{N_{j_s},\ldots, N_{j_1}\}$ in $(*)$ so that we get an exceptional sequence
$$(M_{m_q}, \ldots, M_{m_1}, N_{j_s},\ldots,N_{j_1}, V, W, M_{l_r}', \ldots, M_{l_1}').$$
Now, it follows from the definition of the $V_i, W_i$ that we get an exceptional sequence
$$(M_{m_q}, \ldots, M_{m_1}, V_{n-2},W_{n-2}, N_{j_s},\ldots,N_{j_1}, M_{l_1}', \ldots, M_{l_r}').$$
We claim that for $1 \le u \le q$, we have that $\C_{m_u}=\C(M_{m_u},V_{m_u},W_{m_u})$ is equivalent to $\C(M_{m_u},V_{n-2},W_{n-2})$. Fix such a $u$. Note that there is an exceptional sequence of the form $$(N_{j_s}, \ldots, N_{j_p}, M_{m_u},V_{m_u},W_{m_u}).$$ Now, it follows from Lemma \ref{Lemma5}(2) that $M_{m_u}$ lies in $^\perp N_{j_i}$ for all $1 \le i \le p$. By reflecting, we get the exceptional sequence
$$(M_{m_u},V_{n-2},W_{n-2}, N_{j_s}, \ldots, N_{j_p}).$$
It follows that $\C_{m_u}=\C(M_{m_u},V_{m_u},W_{m_u})$ is equivalent to $\C(M_{m_u},V_{n-2},W_{n-2})$. This proves our claim. Let $E = (M_{m_q}, \ldots, M_{m_1}, V_{n-2},W_{n-2})$. By Lemma \ref{Lemma4} and our claim, $\mathcal{R}(Q,\delta)=\C(E)$
is tame connected. Since $V_{n-2}, W_{n-2}$ lie in it, $\bar \delta$ is the (unique) isotropic Schur root of $\mathcal{R}(Q,\delta)$. It follows from Lemma \ref{Lemma5} that any $\sigma_\delta$-stable representation not isomorphic to any $M_i$ for $1 \le i \le n-2$ will have to be (quasi-simple) in $\C(E)$. Now, we need to show that all quasi-simple objects of $\C(E)$ are $\sigma_\delta$-stable. Assume the contrary. Let $f$ be the dimension vector of a quasi-simple object in $\C(E)$ that is not $\sigma_\delta$-stable, but $\sigma_\delta$-semistable. It follows from our previous observations that $f$ has to be a positive integral combination of the $\sigma_\delta$-stable dimension vectors in $\C(E)$. It follow from \cite[Cor. 6.3]{IPT} that this is not possible to have such a decomposition.
Therefore, we have a complete list of the simple objects in $\A(\delta)$. These are given by the disjoint union
$$\{M_i \mid i \not \in I\} \cup \{\text{quasi-simple objects in} \; \mathcal{R}(Q,\delta)\}.$$
Observe that, in $C_\mathbb{R}(\sigma_\delta)$, a dimension vector $d$ can be uniquely written as
$$d = d_1 + \sum_{i \not \in I} \lambda_i f_i$$
where $d_1$ is a dimension vector in $\C(E)$ and $f_i = d_{M_{m_i}}$ for $i \not \in I$. This decomposition is unique. This implies the unicity of $\mathcal{R}(Q,\delta)$ and statement $(3)$.
\end{proof}

\begin{Cor} Let $Q$ be an acyclic connected quiver and $\delta$ an isotropic Schur root. Then
${\rm SI}(Q, \delta)$ is a polynomial ring or a hypersurface. More precisely, it is a hypersurface (and not a polynomial ring) if and only if $\mathcal{R}(Q,\delta)$ has quiver of type $\widetilde{\mathbb{D}_n}$ with $n \ge 4$, $\widetilde{\mathbb{E}_6}$, $\widetilde{\mathbb{E}_7}$ or $\widetilde{\mathbb{E}_8}$.
\end{Cor}

\begin{proof} In \cite{Skowronskiweyman}, it was proven that the ring of semi-invariant of an isotropic Schur root of a tame quiver is a polynomial ring or a hypersurface, where the second situation occurs precisely when the quiver is of type $\widetilde{\mathbb{D}_n}$ with $n \ge 4$, $\widetilde{\mathbb{E}_6}$, $\widetilde{\mathbb{E}_7}$ or $\widetilde{\mathbb{E}_8}$. Our result follows from this and Theorem \ref{MainTheo}.
\end{proof}

\begin{Exam} Consider the quiver $Q$ given by
$$\xymatrixrowsep{7pt}\xymatrix{& 2 \ar[dl] & \\ 1 & & 4 \ar[ul] \ar[ll] \ar[dl]\\ & 3 \ar[ul]  &}$$
Consider the exceptional sequence $(P_2, S_1, I_3, S_3)$ where $P_2$ is the projective representation at vertex $2$, $I_3$ is the injective representation at vertex $3$ and $S_1, S_3$ are the simple representations at vertices $1,3$, respectively. Reflecting $S_1, I_3$ to the left of $P_2$, we get an exceptional sequence whose dimension vectors are as follows.
$$((0,1,0,0), (3, 3, 1, 1), (1,1,0,0), (0,0,1,0)).$$
Then, using a sequence of reflections, we get the following exceptional sequences, where we put the corresponding dimension vectors.
$$((0,1,0,0), (3, 3, 1, 1), (0,0,1,0), (1,1,1,0))$$
$$((0,1,0,0), (0,0,1,0), (3, 3, 3, 1), (1,1,1,0))$$
$$((0,0,1,0), (0,1,0,0), (3, 3, 3, 1), (1,1,1,0))$$
$$((0,0,1,0), (0,1,0,0), (8,8,8,3), (3, 3, 3, 1))$$
$$((0,0,1,0), (8,3,8,3), (0,1,0,0), (3, 3, 3, 1))$$
$$((8,3,3,3), (0,0,1,0), (0,1,0,0), (3, 3, 3, 1)).$$
Observe that $$\langle (3,3,3,1),(0,1,0,0) \rangle =2$$ and $$\delta = (3,3,3,1) - (0,1,0,0) = (3,2,3,1)$$ is an isotropic Schur root.
The Coxeter matrix $\tau$ is
$$\tau = \left(%
\begin{array}{cccc}
  -1  & 1 & 1 & 1 \\
  -1 & 0 & 1 & 2 \\
  -1 & 1 & 0 & 2 \\
  -3 & 2 & 2 & 4 \\
\end{array}%
\right).$$
This matrix has eigenvalues $\lambda = 5/2 +\sqrt{21}/2$, $\lambda^{-1} = 5/2-\sqrt{21}/2$ and $-1$ with (algebraic and geometric) multiplicity $2$. The eigenvector corresponding to $\lambda$ is $$v_1 = (10, 9+\sqrt{21}, 9+\sqrt{21}, 17 + \sqrt{189})$$ and the one corresponding to $\lambda^{-1}$ is $$v_2 = (10, 9-\sqrt{21}, 9-\sqrt{21}, 17 - \sqrt{189}).$$
Now, $\langle v_2, (8,3,3,3) \rangle = -197 + 10\sqrt{21} + 11\sqrt{189} > 0$ and $\langle v_2, (0,0,1,0) \rangle = -8-\sqrt{21}+\sqrt{189} > 0$.
Similarly, both $\langle (8,3,3,3), v_1 \rangle$ and $\langle (0,0,1,0), v_1 \rangle$ are positive.
Therefore, the exceptional objects with dimension vectors $(8,3,3,3), (0,0,1,0)$ are regular by the theorem at page 240 of \cite{Pena}. It follows from Proposition \ref{Prop4Vertices} that $\delta$ is not of smaller type. It also follows from the same proposition that there is a unique exceptional sequence $(M_2, M_1)$ of length $2$ of $\sigma$-stable objects. Let $M_1' = S_3$ and $M_2'$ be the exceptional representation with dimension vector $(8,3,3,3)$. Since $M_1', M_2'$ lie in $\C(M_2,M_1)$ by Lemma \ref{Lemma0}, we see that $\C(M_2',M_1') \subseteq \C(M_2,M_1)$ and thus, we have equality. This means that $M_2' = M_2$, $M_1' = M_1$. Since $\langle \delta, (0,0,1,0) \rangle = 2 >0$, we get $\delta_1 = \delta - 2(0,0,1,0) = (3,2,1,1)$. Now, $\langle \delta_1, (8,3,3,3) \rangle = -2$. Therefore, $\bar \delta = \delta_1 = (3,2,1,1)$. In this example, the cone of $\sigma$-semistable dimension vectors is as follows (where only an affine slice of that cone is shown).
\end{Exam}

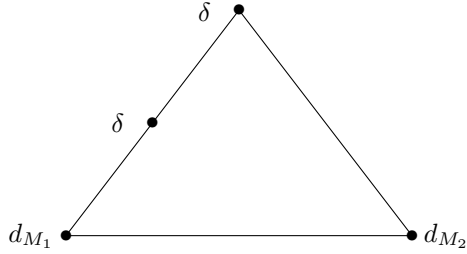
\begin{figure}[h]
  \centering
  \begin{tikzpicture}[xscale=2.30,yscale=1.5]

    \draw (-1,-1) -- (1,-1);

    \node at (-1,-1) {$\bullet$};
\node at (-1.2,-1) {$d_{M_1}$};
    \node at (1,-1) {$\bullet$};
\node at (1.2,-1) {$d_{M_2}$};
\node at (-0.5,0) {$\bullet$};
\node at (-0.7,0) {$\delta$};
\node at (0,1) {$\bullet$};
\node at (-0.2,1) {$\bar \delta$};

    \draw (-1.00,-1) -- (0,1);
\draw (0,1) -- (1,-1);

  \end{tikzpicture}
\caption{The cone of $\sigma$-semistable dimension vectors for $\delta = (3,2,3,1)$}
\label{fig:1}
\end{figure}

The following is an easy observation. The reader is referred to \cite{FZ} for the notion of cluster algebra and to \cite{Fei} for results in similar directions.

\begin{Cor} If ${\rm SI}(Q,\delta)$ is not a polynomial ring, then it has a cluster algebra structure of type $\mathbb{A}_1$. There are two cluster variables which are all $\Gamma$-homogeneous, and the coefficients are built from $n-1$ frozen variables, which are also $\Gamma$-homogeneous, where $\Gamma$ is the set of all multiplicative characters of ${\rm GL}_\delta(k)$.
\end{Cor}

\begin{proof} From Theorem \ref{MainTheo}, it is enough to prove this for $\rep(Q) = \mathcal{R}(Q,\delta)$, that is, we may assume that $Q$ is tame connected. Suppose that ${\rm SI}(Q,\delta)$ is not a polynomial ring. Then $Q$ is of type $\widetilde{\mathbb{D}_n}$ with $n \ge 4$, $\widetilde{\mathbb{E}_6}$, $\widetilde{\mathbb{E}_7}$ or $\widetilde{\mathbb{E}_8}$. In particular, it is well known in these cases that there are exactly three non-homogeneous tubes $T_1, T_2, T_3$ in the Auslander-Reiten quiver of $\mathcal{R}(Q,\delta)$. One, say $T_1$, has rank $2$. Let $M,N$ be the non-isomorphic exceptional quasi-simple objects in $T_1$. Then, let $E_1, \ldots, E_r$ be the non-isomorphic quasi-simple objects of $T_2$ and let $E_1', \ldots, E_t'$ be the non-isomorphic quasi-simple objects of $T_3$. Now, the hypersurface equation can be written as
$$(*) \qquad C^MC^N = C^{E_1}\cdots C^{E_r} + C^{E_1'}\cdots C^{E_t'}.$$
Consider the indeterminates $x,y_1, \ldots, y_r, z_1, \ldots, z_t$. We define a cluster algebra $A$ as follows. We start with the initial seed $\{x,y_1, \ldots, y_r, z_1, \ldots, z_t\}$ where $y_1, \ldots, y_r$ and $z_1, \ldots, z_t$ are declared to be frozen variables. The exchange relation is
$xx' = \prod_{i=1}^ry_i + \prod_{j=1}^tz_j$
which clearly produces exactly two cluster variables $x,x'$. The cluster algebra is the the $\Z$-subalgebra of $\mathbb{Q}(x,y_1, \ldots, y_r, z_1, \ldots, z_t)$ generated by $x, x'$ and $y_1, \ldots, y_r, z_1, \ldots, z_t$. This algebra is clearly isomorphic to ${\rm SI}(Q,\delta)$.
\end{proof}

An interesting problem would be to find all acyclic quivers $Q$ and dimension vectors $d$ such that SI$(Q,d)$ has a cluster algebra structure whose variables (frozen or not) are all $\Gamma$-homogeneous.

\section{Construction of all isotropic Schur roots}

In this section, we show that all of the isotropic Schur roots of $\rep(Q)$ come from isotropic Schur roots of a tame full subquiver of $Q$ by applying special reflections. We make this precise by defining an action of the braid group $B_{n-1}$ on $n-1$ strands on a special type of exceptional sequences that will encode all we need to study isotropic Schur roots. We start with the definition of these sequences.

\begin{Defn}
Let $E=(X_1, \ldots, X_n)$ be a full exceptional sequence. We say that $E$ is of \emph{isotropic type} if there exists $1 \le i \le n-1$ such that $\C(X_{i}, X_{i+1})$ is tame.
The integer $i$ is called the \emph{isotropic position} of $E$ and the \emph{root type} of $E$, denoted $\delta_E$, is the isotropic Schur root in $\C(X_{i}, X_{i+1})$.
\end{Defn}

We denote by $\mathcal{E}$ the set of all full exceptional sequences of isotropic type, up to isomorphism. Not all elements of the braid group $B_{n}$ act on $\mathcal{E}$. We rather consider the group $B_{n-1}$ and show that it acts on $\mathcal{E}$. Let us denote the standard generators of $B_{n-1}$ by $\gamma_1, \ldots, \gamma_{n-2}$. Let $E=(X_1, \ldots, X_n) \in \mathcal{E}$ with isotropic position $r$. Let $1 \le i \le n-2$. If $i<r-1$, then $\gamma_iE := \sigma_iE$. If $i>r$, then $\gamma_{i}E := \sigma_{i+1}E$. Assume that $i=r$ with $r<n-1$. We can reflect $X_{r+2}$ to the left of $X_r, X_{r+1}$ to get the exceptional sequence:
$$E'=(X_1, \ldots,X_{r-1}, L_{X_{r}}(L_{X_{r+1}}(X_{r+2})), X_r, X_{r+1},  X_{r+3}, \ldots, X_n).$$
and this is an exceptional sequence of isotropic type with isotropic position $r+1$. We define $\gamma_rE:=E'$. If $r > 1$ and $i=r-1$, then we can reflect both $X_r, X_{r+1}$ to the left of $X_{r-1}$ as follows:
$$E'' = (X_1, \ldots,X_{r-2}, L_{X_{r-1}}(X_r),  L_{X_{r-1}}(X_{r+1}), X_{r-1},  X_{r+2}, \ldots, X_n)$$
and clearly, the subcategory $\C(L_{X_{i-1}}(X_i),  L_{X_{i-1}}(X_{i+1}))$ generates a tame subcategory of rank $2$. Therefore, $E''$ is an exceptional sequence of isotropic type with isotropic position $r-1$ and its root type is the unique isotropic Schur root in $\C(L_{X_{i-1}}(X_i),  L_{X_{i-1}}(X_{i+1}))$, which is $\delta_{E'} = \delta_E - \langle \delta_E, d_{X_{r-1}}\rangle d_{X_{r-1}}$, by Lemma \ref{LemmaIsoPreinj}. We define $\gamma_{r-1}E=E''$.
Similarly, we can define the action of $\gamma_{i}^{-1}$ on $E$ for $1 \le i \le n-2$. The following is easy to check.

\begin{Prop} The group $B_{n-1}$ acts on exceptional sequences of isotropic type, with the action defined above.
\end{Prop}

\begin{Defn} A sequence $E=(X_1, \ldots, X_{n-1}, X_n)$ in $\mathcal{E}$ is of \emph{tame type} if it has isotropic position $n-1$, and there is $0 \le s \le n-2$ such that $X_1, \ldots, X_s$ are projective in $\rep(Q)$ and $\C(X_{s+1}, \ldots, X_{n-2}, X_{n-1}, X_n)$
is tame connected. By convention, $s=0$ means that $\rep(Q)$ is already tame connected.
\end{Defn}

Observe that if $E \in \mathcal{E}$ is of tame type, then the isotropic Schur root $\delta_E$ is the unique isotropic Schur root of the tame subcategory $\C(X_{s+1}, \ldots, X_{n-2}, X_{n-1}, X_n)$ and is an isotropic Schur root coming from a tame full subquiver of $Q$. In particular, there are finitely many roots $\delta_E$ where $E \in \mathcal{E}$ is of tame type.

\begin{Exam} Consider a quiver of rank $n=4$ and an exceptional sequence $E=(X,U,V,Y)$ of isotropic type with isotropic position $2$. The root type is the isotropic root $\delta_E$ in $\C(U,V)$.
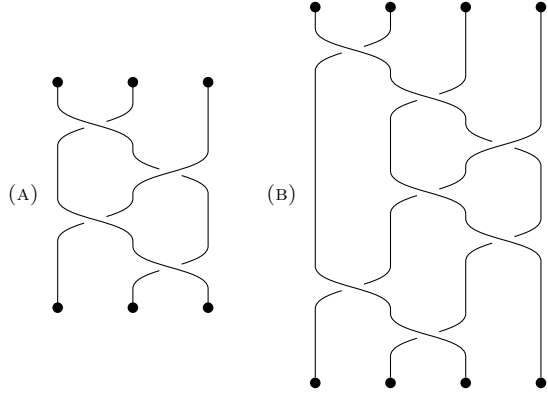
\begin{figure}[h]
\sidesubfloat[]{%
\begin{tikzpicture}
\braid[height = 18pt] a_1a_2^{-1}a_1a_2;
 \fill [black] (1,0) circle (2pt);
\fill [black] (2,0) circle (2pt);
\fill [black] (3,0) circle (2pt);
\fill [black] (1,-3) circle (2pt);
\fill [black] (2,-3) circle (2pt);
\fill [black] (3,-3) circle (2pt);
\end{tikzpicture}}\qquad%
\sidesubfloat[]{%
\begin{tikzpicture}
\braid[height = 18pt] a_1a_2a_3^{-1}a_2a_3a_1a_2;
 \fill [black] (1,0) circle (2pt);
\fill [black] (2,0) circle (2pt);
\fill [black] (3,0) circle (2pt);
\fill [black] (4,0) circle (2pt);
\fill [black] (1,-5) circle (2pt);
\fill [black] (2,-5) circle (2pt);
\fill [black] (3,-5) circle (2pt);
\fill [black] (4,-5) circle (2pt);
\end{tikzpicture}}
\caption{Correspondence between $B_3$ and some braids of $B_4$}
\label{fig:2}
\end{figure}

\noindent The first braid (A) in Figure \ref{fig:2} corresponds to the element $g=\gamma_2^{-1}\gamma_1^{-1}\gamma_2\gamma_1^{-1}$ of $B_3$ while the second braid (B) corresponds to the element $h=\sigma_2^{-1}\sigma_1^{-1}\sigma_3^{-1}\sigma_2^{-1}\sigma_3\sigma_2^{-1}\sigma_1^{-1}$ of $B_4$. Notice that $gE = hE$. Notice also that the braid in (A) is obtained from the braid in (B) by identifying the two strands starting at the positions of $U,V$, that is, the second and third strands.
\end{Exam}

Our aim in this section is to prove that any $E \in \mathcal{E}$ lies in the $B_{n-1}$-orbit of an exceptional sequence of tame type. In the next lemmas, we will consider exceptional sequences in the bounded derived category $D^b(\rep(Q))$ of $\rep(Q)$. Recall that an object $X$ in $D^b(\rep(Q))$ is \emph{exceptional} if $\Hom(X,X[i])=0$ for all non-zero $i$ (and then, $\Hom(X,X)$ has to be one dimensional). Equivalently, an exceptional object in $D^b(\rep(Q))$ is isomorphic to the shift of an exceptional representation. A sequence $(X_1, \ldots, X_r)$ of objects in $D^b(\rep(Q))$ is \emph{exceptional} if every $X_i$ is exceptional and, for $i < j$, we have $\Hom_{D^b(\rep(Q))}(X_i, X_j[t])=0$ for all $t \in \Z$. For such a sequence, one can consider the smallest full additive subcategory $\mathcal{D}(X_1, \ldots, X_r)$ of $D^b(\rep(Q))$ containing $X_1, \ldots, X_r$ and that is closed under direct sums, direct summands, taking the cone of a morphism and the shift of an object. One can also consider the exceptional sequence $(X_1', \ldots, X_r')$ in $\rep(Q)$ such that $X_i'$ is the unique shift of $X_i$ lying in $\rep(Q)$. The indecomposable objects in $\mathcal{D}(X_1, \ldots, X_r)$ are just the shifts of the indecomposable objects in $\C(X_1', \ldots, X_r')$.

\medskip

In what follows, the Auslander-Reiten translate in $D^b(\rep(Q))$ is denoted by $\tau_{D}$ while the Auslander-Reiten translate in $\rep(Q)$ is simply denoted $\tau$. Recall that if $X$ is a non-projective indecomposable representation, then $\tau_DX = \tau X$ and, if $X = P_x$ with $x \in Q_0$, then $\tau_D X = I_x[-1]$. When $d$ is a dimension vector, we denote by $\tau d$ the product of the Coxeter matrix with $d$. In particular, if $X$ is a non-projective indecomposable representation, then $\tau d_X = d_{\tau X}$ and, if $X=P_x$ with $x \in Q_0$, then $\tau d_{X} = -d_{I_x}$.  We start our investigation with the following lemma that is crucial for the proof of the main result of this section.

\begin{Lemma} \label{Lemma8} Let $(X_1, \ldots, X_{n})$ be an exceptional sequence with $\C(X_{r+1}, \ldots, X_{n})$ tame and assume that $X_1, \ldots, X_r$ are the simple objects in $\C(X_1, \ldots, X_r)$. Let $X \in \C(X_1, \ldots, X_r)$ be the injective object with socle $X_1$. If $X$ is projective in $\C(X,X_{r+1}, \ldots, X_{n})$, then $X_1$ is projective in $\rep(Q)$ and in particular, an isotropic Schur root of $\C(X_{r+1}, \ldots, X_{n} )$ is not sincere.
\end{Lemma}

\begin{proof} Assume that $X$ is projective in $\C(X,X_{r+1}, \ldots, X_{n})$.  Set $d_i = d_{X_i}$ for $1 \le i \le n$. Consider the linear form $f$ given by $f(x)=\langle d_1, x \rangle$. Then $f$ vanishes on $d_2, \ldots, d_{n}$ and $f(d_1)>0$. Assume to the contrary that $X_1$ is not projective in $\rep(Q)$. Observe that $f(x)=\langle d_1, x \rangle = -\langle x, \tau d_1 \rangle$. Since $\tau X_1$ is exceptional, $\langle \tau d_1, \tau d_1 \rangle = 1$ and hence $f(\tau d_1) < 0$. Now, reflect $X_1$ to the right of $X_2, \ldots X_r$, so that we get an exceptional sequence $(X_2, \ldots, X_r, Y)$ where $Y$ is in the cone spanned by $d_1, \ldots, d_r$. Clearly, $X_1$ is simple projective in $\C(X_1, \ldots, X_r)$ and hence, $Y=X$ is the injective hull of $X_1$ in $\C(X_1, \ldots, X_r)$. Set $\C:=\C(X, X_{r+1}, \ldots, X_{n})$. We know that $X$ is projective in $\C$. Reflecting $X$ to the right of $X_{r+1}, \ldots, X_n$ will give the exceptional representation $\tau X_1$. Therefore, $d:=\tau d_1 = -\tau_\C d_X$ where $\tau_\C$ denotes the Coxeter transformation in $\C$. Take the linear form $g$ in the Grothendieck group of $\C$ given by $g(x)=\langle d_{X}, x \rangle$. Then $g$ vanishes on $d_{r+1}, \ldots, d_{n}$ and $g(d_X)>0$. The form $f|_\C$ has the same property since $d_X$ is a non-negative linear combination of $d_1, \ldots, d_r$ with the coefficient of $d_1$ positive. Thus, $g = f|_\C$ up to a positive scalar. Therefore, $g(d)<0$, which means that $X,\tau X_1$ lie on opposite sides of the hyperplane $g(x)=0$ in $\C$. This contradicts that $X$ is projective in $\C$.
\end{proof}

\begin{Lemma} \label{MinimalRoot}
Let $(X_1, \ldots, X_r)$ be an exceptional sequence and assume that $\C_1:=\C(X_2, \ldots, X_r)$ is tame with an isotropic Schur root $\gamma$ while $\C_2:=\C(X_1, X_2, \ldots, X_r)$ is wild. Then there is a unique minimal isotropic Schur root in the $\tau$-orbit of $\gamma$.
\end{Lemma}

\begin{proof} We may assume that $\C_2 = \rep(Q)$ for an acyclic quiver $Q$. Since $\C_2$ is wild and $\C_1$ is tame, we know that $X_1$ is preprojective or preinjective in $\C_2$. Hence, there is some $r \in \Z$ such that $\tau_D^r X_1$ is projective or the shift of a projective. This means that $\tau^r\gamma$ is not sincere. Let $Y = \tau_D^r X_1$ if $\tau_D^r X_1$ is a representation or $Y = \tau_D^r X_1[-1]$ if $\tau_D^r X_1$ is the shift of a projective representation. Observe that $Y^\perp \subseteq \rep(Q)$ is also of tame representation type, where the quivers of $X_1^\perp$ and $Y^\perp$ only differ by a change of orientation; see for instance \cite[Prop. 2.1]{Happel}. Therefore, $\tau^r\gamma$ is an isotropic Schur root of a tame full subquiver of $Q$. Let $s \in \Z$ with $s \ne r$. Consider $Z$ the unique shift of $\tau_D^s X_1$ which is a representation. Since the simples in $Y^\perp$ are simples in $\rep(Q)$, and since there is a simple of $Z^\perp \subseteq \rep(Q)$ that is not simple in $\rep(Q)$, we see that the isotropic Schur root $\tau^r\gamma$ has smaller length than $\tau^s\gamma$. This also proves unicity since only one object in the $\tau$-orbit of $X_1$ is projective or a shift of a projective.
\end{proof}

\begin{Lemma} \label{Lemma6}Let $E=(X_1, \ldots, X_{n-2}, U, V)$ be in $\mathcal{E}$ with isotropic position $n-1$. Let $E' = \gamma_{1}\cdots\gamma_{n-3}\gamma_{n-2}E = (U', V', X_1, \ldots, X_{n-2})$.
Then $\tau^{-1} \delta_E = \delta_{E'}$.
\end{Lemma}

\begin{proof}
If $V$ is not injective, then we have the exceptional sequence $$(\tau^{-1}V, X_1, \ldots, X_{n-2},U)$$ in $\rep(Q)$. Otherwise, we have the exceptional sequence $$(\tau^{-1}V[-1], X_1, \ldots, X_{n-2}, U)$$ in $\rep(Q)$. Let us write $\tau^{-1}V[0,-1]$ to indicate that we either take the shift $[0]$ or $[-1]$ for $\tau^{-1}V$. Then, we get an exceptional sequence $$(\tau^{-1}U[0,1], \tau^{-1}V[0,1], X_1, \ldots, X_{n-2}).$$ The categories $\C(U', V')$ and $\C(\tau^{-1}U[0,1], \tau^{-1}V[0,1])$ are equal in $\rep(Q)$. Therefore, they have the same isotropic Schur root. The isotropic Schur root of $$\C(\tau^{-1}U[0,1], \tau^{-1}V[0,1])$$ is clearly $\tau^{-1}\delta_E$.
\end{proof}

Of course, we have the dual version of the above lemma as follows.

\begin{Lemma} \label{Lemma7}
Let $E=(U,V, X_1, \ldots, X_{n-2})$ be in $\mathcal{E}$ with isotropic position $1$. Let $E' = \gamma_{n-2}^{-1}\gamma_{n-3}^{-1}\cdots\gamma_{1}^{-1}E = (X_1, \ldots, X_{n-2}, U', V')$.
Then $\tau \delta_E = \delta_{E'}$.
\end{Lemma}

We are now ready for the main result of this section.

\begin{Theo}
Let $\delta$ be an isotropic Schur root. Then there is $E \in \mathcal{E}$ of tame type and $g \in B_{n-1}$ such that $gE$ has root type $\delta$.
\end{Theo}

\begin{proof}
It follows from Proposition \ref{ExceptionalSequenceIsotropic} that there is an exceptional sequence $F=(M_1, \ldots, M_{n-2}, X,Y)$ in $\mathcal{E}$ of isotropic position $n-1$ and of root type $\delta$. Assume that $G\in \mathcal{E}$ is in the orbit of $E$ and the root type of $G$ is minimal, that is, has minimal length as a root in $\rep(Q)$. We may assume that the isotropic position of $G$ is $n-1$.
Therefore, we may assume that $G$ is of the form
$$(Y_1, \ldots, Y_{n-2}, U, V).$$
Assume first that there is an object $W$ in $\C(Y_1, \ldots, Y_{n-2})$ that is not projective in $\C(W,U,V)$. We can apply a sequence of reflections to the subsequence $(Y_1, \ldots, Y_{n-2})$ to get an exceptional sequence
$$H=(Y_1', \ldots, Y_{n-3}', W, U,V)$$
in $\mathcal{E}$.
Now, applying $\gamma_{n-2}^2$ to $H$ and using Lemma \ref{Lemma6}, we get the sequence $$(Y_1', \ldots, Y_{n-3}', W', U', V')$$ in $\mathcal{E}$ whose root type is the inverse Auslander-Reiten translate of $\delta$ in $\C(W,U,V)$.
Similarly, applying $(\gamma_{n-2}^{-1})^2$ to $H$, we get the sequence
$$(Y_1', \ldots, Y_{n-3}', W'', U'', V'')$$
in $\mathcal{E}$ whose root type is the Auslander-Reiten translate of $\delta$ in $\C(W,U,V)$.
We can iterate this to get a smaller root by Lemma \ref{MinimalRoot}, provided $\C(W,U,V)$ is wild. Therefore, whenever there is an object $W$ which is not projective in $\C(W,U,V)$, then $\C(W,U,V)$ is of tame type (and hence connected). Suppose, by induction, that we have an exceptional sequence $J=(W_{r+1},\ldots,W_{n-2},U,V)$ such that $\C(J)$ is tame connected and $J$ has maximal length with respect to this property. If $r=0$, then $Q$ is a tame connected quiver and there is nothing to prove. Complete this to get a full exceptional sequence
$$(Z_1, \ldots, Z_r, W_{r+1}, \ldots, W_{n-2}, U,V).$$
If there is $W \in ^\perp \C(J) = \C(Z_1, \ldots, Z_r)$ such that $W$ is not projective in $\C(W,U,V)$ then, $\C(W,U,V)$ is tame connected. As in the proof of Lemma \ref{Lemma4}, we get that $\C(W,W_{r+1}, \ldots, W_{n-2},U,V)$ is tame connected, contradicting the maximality of $J$. Therefore, any object $Z$ in $\C(Z_1, \ldots, Z_r)$ is such that $Z$ is projective in $\C(Z,W_{r+1}, \ldots, W_{n-2},U,V)$. We may apply a sequence of reflections and assume that all of $Z_1, \ldots, Z_r$ are simple in $\C(Z_1, \ldots, Z_r)$. It follows from Lemma \ref{Lemma8} that the injective hull of $Z_1$ in $\C(Z_1, \ldots, Z_r)$ is projective in $\rep(Q)$. Then the proof goes by induction.
\end{proof}

Here is another way to interpret this result. Start with an isotropic Schur root $\delta_0$ of a tame full subquiver $Q'$
of $Q$ and consider an exceptional sequence $(U_0,V_0)$ of length $2$ in $\rep(Q') \subset \rep(Q)$
such that $\delta_0 = d_{U_0} + d_{V_0}$. Consider an exceptional object $X_0$ such that $(X_0,U_0,V_0)$ is an exceptional sequence of length three (which generates a thick subcategory $\C_0$ of $\rep(Q)$). Then we can transform it into another exceptional sequence $(X_0', U_1, V_1)$ with an isotropic Schur root $\delta_1 = d_{U_1}+d_{V_1}$ such that $\delta_1$ is a power $\tau_{\C_0}^{r_0} \delta_0$ where $\tau_{\C_0}$ denotes the Coxeter matrix for $\C_0$. Now, for $i \ge 1$, consider an exceptional object $X_i$ such that $(X_i,U_i, V_i)$ is an exceptional sequence. Take a power $\delta_{i+1}=\tau_{\C_i}^{r_i} \delta_i$ where $\C_i$ is the thick subcategory of $\rep(Q)$ generated by $X_i, U_i, V_i$ and $\tau_{\C_i}$ denotes the Coxeter matrix for $\C_i$. All the roots $\delta_i$ constructed this way are isotropic Schur roots. Moreover, all isotropic Schur roots of $\rep(Q)$ can be obtained in this way. There are clearly only finitely many starting roots $\delta_0$, but the choices of the $r_i$ and $X_i$ yield, in general, infinitely many possible isotropic Schur roots. As observed in \cite{Pa}, when $Q$ is wild connected with more than $3$ vertices, there are infinitely many $\tau$-orbit of isotropic Schur roots (provided there is at least one isotropic Schur root). An interesting question would be to describe the minimal root types of the orbits of $\mathcal{E}$ under $B_{n-1}$. It is not hard to check that when $n=3$, these minimal root types correspond exactly to the tame full subquivers of $Q$. We do not know if this holds in general.

\begin{Conj}
Let $E_1, E_2 \in \mathcal{E}$. Assume that there are $g_1, g_2 \in B_{n-1}$ with $g_1E_1, g_2E_2$ of tame type but with different root types. Then $E_1, E_2$ lie in distinct orbits under $B_{n-1}$.
\end{Conj}

\bigskip

\noindent{\emph{Acknowledgment}.} The authors are thankful to Hugh Thomas for suggesting the decomposition in Proposition \ref{PropVectSpacesDecomp}. The second named author was supported by NSF grant DMS-1400740.

\end{document}